\documentclass[preprint,11pt]{article}

\usepackage{lineno,hyperref}

\usepackage{graphicx, amssymb, amsmath}
\usepackage{indentfirst}
\usepackage{mathrsfs}
\usepackage{amsfonts}
\usepackage{color}
\usepackage{upgreek}
\usepackage{subfigure}
\usepackage[misc]{ifsym}
\usepackage{hyperref}
\usepackage{float}
\usepackage{xcolor}

\usepackage{epstopdf}
\usepackage{booktabs}
\usepackage{multirow}
\usepackage{makecell}

\usepackage{geometry}
\usepackage{algorithmic}   
\usepackage[ruled,linesnumbered]{algorithm2e}

\usepackage[noblocks]{authblk} 

\newtheorem{example}{\emph{Example}}[section]

\newcommand{\be}{\begin{equation}}
\newcommand{\ee}{\end{equation}}
\newcommand{\ba}{\begin{array}}
\newcommand{\ea}{\end{array}}
\newcommand{\beas}{\begin{eqnarray*}}
\newcommand{\eeas}{\end{eqnarray*}}
\newcommand{\bea}{\begin{eqnarray}}
\newcommand{\eea}{\end{eqnarray}}

\begin{document}

\title{ \Large An inverse averaging finite element method for solving the size-modified Poisson-Nernst-Planck equations in ion channel simulations}

\author[1]{Ruigang Shen}
\author[1,2]{Qianru Zhang}
\author[1,\thanks{Corresponding author. \newline
	 Email addresses: shenruigang@lsec.cc.ac.cn (S.~Shen),  qrzhang@lsec.cc.ac.cn (Q.~Zhang),
     bzlu@lsec.cc.ac.cn (B.~Lu)} ] { Benzhuo Lu }

\affil[1]{\footnotesize ~LSEC, Institute of Computational Mathematics and Scientific/Engineering Computing, Academy of Mathematics and Systems Science, Chinese Academy of Sciences, Beijing, 100190, P.R. China.}

\affil[2]{ School of Mathematical Sciences, University of Chinese Academy of Sciences, Beijing 100049, P.R. China.}


\maketitle


\noindent
\textbf{Abstract}:
  In this work, we introduce an inverse averaging finite element method (IAFEM) for solving the size-modified Poisson-Nernst-Planck (SMPNP) equations. Comparing with the classical Poisson-Nernst-Planck (PNP) equations, the SMPNP equations add a nonlinear term to each of the Nernst-Planck (NP) fluxes to describe the steric repulsion which can treat multiple nonuniform particle sizes in simulations. Since the new terms include sums and gradients of ion concentrations, the nonlinear coupling of SMPNP equations is much stronger than that of PNP equations.
  By introducing a generalized Slotboom transform, each of the size-modified NP equation is transformed into a self-adjoint equation with exponentially behaved coefficient, which has similar simple form to the standard NP equation with the Slotboom transformation. This treatment enables employing our recently developed inverse averaging technique to deal with the exponential coefficients of the reformulated formulations, featured with advantages of numerical stability and flux conservation especially in strong nonlinear and convection-dominated cases.
  Comparing with previous stabilization methods, the IAFEM proposed in this paper can still possess the numerical stability when dealing with convection-dominated problems. And it is more concise and easier to be numerically implemented. Numerical experiments about a model problem with analytic solutions are presented to verify the accuracy and order of IAFEM for SMPNP equations. Studies about the size-effects of a sphere model and an ion channel system are presented to show that our IAFEM is more effective and robust than the traditional finite element method (FEM) when solving SMPNP equations in simulations of biological systems. \\



\section{Introduction}
\label{sec1}

\noindent
``The effects of finite particle size on electrostatics, density profiles, and diffusion have been a long-existing topic in the study of ionic solution."\cite{Lu2011smpnp} As a continuous electrodiffusion model, the classical Poisson-Nernst-Planck (PNP) equations play an important role in the electrodiffusion reaction process and have been widely used to describe the electrodiffusion of ions and charge transport in applications including the solvated biomolecular system \cite{Lu2010pnp,Lu2007Electrodiffusion}, semiconductors \cite{Jerome1996,Markowich1986,Selberherr1984}, electrochemical systems \cite{Bazant2009Towards,Ciucci2011Derivation,Marcicki2012Comparison} and ion channels \cite{AE2000Three,Coalson2005Poisson,Singer2008A}. Although the PNP equations have achieved a lot of success in various applications, it still has some limitations due to the neglected steric effects of ions in its mean-field derivation, for example, the PNP model leads to unphysical crowding of ions near charged surfaces and incorrect dynamics of ion transport, and the difference between two cations with the same charge cannot be distinguished when simulating the concentration distribution of ions. To incorporate the effects of finite particle sizes in the study of ionic solutions, many improvements are made through introducing exclusion terms from the liquid-state theory or the density functional theory (DFT), e.g. see \cite{Gillespie2002Coupling,Gillespie2003Density,Outhwaite1980Theory,Rosenfeld1997Fundamental} and references therein. In addition, based on the framework of the PNP model, several versions of the modified PNP theory have been developed in the literature to account for steric effects \cite{Horng2012PNP,Kilic2007Steric,Li2009Continuum,Lin2014A,Qiao2016A,Siddiqua2017A}. Among these theories, the Borukhov model \cite{Borukhov1998Steric} attracts people's attention because it captures basic size effects only with a simplified model. The Borukhov model modifies the free energy functional of the ionic system (mean-field approximation) by adding an ideal-gas-like solvent entropy term, which represents the unfavorable energy used to model the over-packing or crowding of the ions and solvent molecules. Thus the steric effects are taken into account in the model. Lu and Zhou by generalizing the Borukhov model get a class of size-modified Poisson-Nernst-Planck (SMPNP) equations via the inclusion of the entropy of solvent molecules in the electrostatic free-energy functional \cite{Lu2011smpnp}. Different from many other works, the SMPNP model is able to treat multiple nonuniform particle sizes in simulations.

 Comparing with the classical PNP equations \cite{Lu2010pnp}, the SMPNP model adds a nonlinear term to each of Nernst-Planck (NP) equations aiming at describing the steric repulsion (see Eq. \eqref{smpnp-np-flux}). Since the new term includes sums and gradients of ion concentrations, the nonlinear coupling of SMPNP equations is much stronger than that of PNP equations. This brings many difficulties for solving SMPNP equations. The NP equations are typical convection-diffusion equations. And the convection dominance will lead to numerical oscillations (e.g. negative ion concentration values). Many stabilization schemes are proposed to avoid non-physical numerical oscillations, e.g. see \cite{Chaudhry2014A,Tu2015Stabilized,Wang2021A}. Tu et al. \cite{Tu2015Stabilized} employed the streamline upwind/Petrov-Galerkin (SUPG) method and the Pseudo Residual-Free Bubble function (PRFB) scheme to enhance the numerical robustness and convergence of the finite element scheme. However, for some macro-molecular systems, e.g., the KcsA ion channel, the SUPG method cannot eliminate all non-physical numerical oscillations or produce convergent numerical solutions \cite{Wang2021A}. By combining the ``upwind" characteristic of the SUPG method and the polishing effect of the interior penalty (IP) method \cite{Burman2004Edge,Douglas1976Interior}, Wang et al.\cite{Wang2021A} proposed a SUPG-IP method to solve PNP equations, which performs better in preserving numerical solution positivity and is much more robust than the standard FEM and the SUPG method when simulating KcsA ion channels. For modified PNP equations with steric effect, based on the ``SUPG" framework, a fast stabilized finite element method is proposed for solving the modified PNP equations with uniform particle sizes \cite{Chaudhry2014A}. However, in the above stabilization schemes, derivations of stabilization terms need a lot of complex interface jump integral calculations, which increases the complexity of their numerical implementations, especially with irregular geometric biological channels. In addition, selecting appropriate stabilization parameters is skillful for different macro-molecule systems. Large stabilization parameters are helpful to numerical convergence, but they will affect numerical accuracy. In addition, the stabilized methods generally have no flux conservation properties.

 Especially, we noticed that the modified PNP equations with ionic steric effects (SPNP) were studied by finite difference methods based on harmonic-mean approximations to the exponential coefficients of the reformulated NP equations in \cite{Ding2019Positivity}. The difference from \cite{Ding2019Positivity} is that the exponential coefficients are approximated based on the inverse averaging of the integral on the element instead of directly on the entire grid node in this paper. Although the finite difference method has been widely used to solve the PNP equations, e.g. see \cite{Flavell2014A,He2016An,Liu2014A} and references therein, the implementation and accuracy of the numerical solution are not so good when it is applied to simulate the actual biomolecular systems with highly irregular surfaces, such as cell membrane, DNA and ion channels. The Finite element method (FEM) has more flexibility and adaptability in irregular regions, which has shown great advantages in solving PNP and modified PNP equations in many actual biomolecular simulations \cite{Liu2017Incorporating,Lu2010pnp,Lu2011smpnp,Song2004Continuum,Xu2014Modeling}.

 In this work, we notice that distributions of ion concentrations cannot be approximated with piecewise polynomials directly, but the size-modified flux densities vary moderately in biological channels. Thus we treat each of the size-modified flux densities as a whole by introducing a set of generalized  Slotboom variables, which eliminate cross-terms in the size-modified NP (SMNP) equations. Then SMNP equations are transformed into self-adjoint second-order elliptic equations with exponentially behaved coefficients. In order to deal with these exponential coefficients, we employ an inverse averaging technique introduced in \cite{ZhangQianru2021A}.
This method possesses good convergence performances when solving
the three-dimensional drift-diffusion (DD) model in semiconductor device simulations and three-dimensional PNP equations in simulating nanopore systems \cite{Zhanng2021An}. It can solve the non-physical spurious oscillation problems caused by the convection domination and guarantee the conservation of computed total currents. Inspired by \cite{ZhangQianru2021A}, we apply the inverse averaging finite element method (IAFEM) for solving SMPNP equations and derive an effective and robust numerical scheme for biomolecular system simulations with the SMPNP model. The main idea is to use the exponential coefficients' harmonic averages to approximate them on every tetrahedral element. And their harmonic averages are calculated on each edge of the tetrahedral element with an inverse averaging technique. This strategy is reasonable because it takes advantage of the moderate variations of size-modified flux densities.

The rest of the paper is organized as follows. In Section \ref{model-intro}, the mathematical model and relevant governing equations are introduced. The singular decomposition of permanent charges for Poisson equations and the reformulation of SMNP equations with the help of generalized Slotboom variables are also presented in this section. In Section \ref{fem-eafe}, the IAFEM are introduced to discretize the SMNP equations. Numerical experiments about a model problem with analytic solutions are presented to verify the accuracy and order of the IAFEM for SMPNP equations. And studies about the size-effects of a sphere model and the an ion channel system are also reported in Section \ref{numerical-experiment}. This paper is ended with Section \ref{sec-conclusion}.

\section{ The mathematical model }
\label{model-intro}

In this section we briefly overview the problem and review the relevant equations.

\subsection{Governing equations }

Let $\Omega \subset \mathbb{R}^d$ $(d=2,3)$ be an open domain. In this work, we consider the SMPNP equations \cite{Lu2011smpnp} by coupling the SMNP equations
\begin{align}
  \frac{\partial c_{i}}{\partial t} &= -\nabla\cdot \boldsymbol{J}_i, \ i = 1, 2,\cdots,K, \ {\rm in} \ \Omega_s, \label{smpnp-np} \\
 \boldsymbol{J}_i &= - D_{i}\Big(\nabla c_{i} + \beta q_ic_{i}\nabla\phi + \frac{k_{i} c_{i} }{1-\sum\limits_{l} a_{l}^{3} c_{l}}\sum\limits_{l} a_{l}^{3} \nabla c_{l} \Big), \label{smpnp-np-flux}
\end{align}
 and the Poisson equation with the internal interface $\Gamma_m = \bar{\Omega}_{s} \cap \bar{\Omega}_{m}$:
 \begin{align}
   &-\nabla\cdot(\epsilon\nabla\phi) = \rho^f + \lambda\sum\limits_{i=1}^K q_i c_i, \ \ {\rm in} \ \Omega=\Omega_s\cup\Omega_m,  \label{smpnp-poi} \\
   &\phi_{m}=\phi_{s}, \quad \epsilon_{m} \frac{\partial \phi_{m}}{\partial \vec{n} }=\epsilon_{s} \frac{\partial \phi_{s}}{\partial  \vec{n}}, \quad x \in \Gamma_m, \label{smpnp-poi-gamma}
 \end{align}
 where $c_{i}(x, t)$ is the concentration of the $i$th ion species carrying charge $q_i = z_i e_c$, $z_i$ is the valence of the $i$th ion species, and $e_c$ is the elementary charge. $\boldsymbol{J}_i$ is the size-modified flux density, in which $D_{i}$ is a spatial-dependent diffusion coefficient, and $\phi$ is the electrostatic potential, $K$ is the number of diffusive ion species considered in the solution system. The constant $\beta=1 /\left(k_{\mathrm{B}} T\right)$ is the inverse Boltzmann energy, where $k_{\mathrm{B}}$ is the Boltzmann constant, and $T$ is the absolute temperature. $\epsilon$ is the piecewise dielectric constant with $\epsilon=\epsilon_{m} \epsilon_{0}$ in $\Omega_{m}$ and $\epsilon=\epsilon_{s} \epsilon_{0}$ in $\Omega_{s}$, where $\varepsilon_{0}$ is the dielectric constant of vacuum, and the typical values of $\epsilon_{m}$ and $\epsilon_{s}$ are 2 and 80, respectively. The characteristic function $\lambda$ is the indicator function of $\Omega_s$, that is
 $\lambda = \left\{
   \begin{array}{l}
             0, \ {\rm in} \ \Omega_m \\
             1, \ {\rm in} \ \Omega_s
   \end{array} \right.$,
which suggests that mobile ions only exist in the solvent region. The permanent (fixed) charge distribution
$$
\rho^{f}(x)=\sum_{j} q_{j} \delta\left(x-x_{j}\right),
$$
which is a sum of singular charges $q_j$ located at $x_j$ inside the biomolecule, and $\delta$ is the Dirac-delta function. The constant $k_{i}=a_{i}^{3} / a_{0}^{3}$, where $a_{i}$ is the effective size of the $i$th ion species, and $a_{0}$ is the solvent molecule size. The size $a_i$ can be arbitrary, and does not need to be larger than the solvent molecule size $a_{0}$.

Comparing with the classical PNP equations (cf. \cite{Lu2010pnp}),  the SMPNP equations add nonlinear terms  $\frac{k_{i} c_{i} }{1-\sum\limits_{l} a_{l}^{3} c_{l}}\sum\limits_{l} a_{l}^{3} \nabla c_{l} $, ($l=1,2,\cdots, K$)  to flux densities in \eqref{smpnp-np-flux} to describe the steric repulsion. If size-effects are not considered, that is, $k_i = 0$ (or $a_i = 0$), SMNP equations \eqref{smpnp-np} directly reduce to classical NP equations as follows
 \begin{align}\label{pnp-np}
  \frac{\partial c_{i}}{\partial t} = \nabla \cdot D_{i}\big(\nabla c_{i} + \beta q_ic_{i}\nabla\phi \big), \ i = 1, 2,\cdots,K, \ {\rm in} \ \Omega_s.
 \end{align}
For brevity, we write the NP equations and SMNP equations as
\begin{align} \label{smpnp-pnp}
  \frac{\partial c_i}{\partial t} &= -\nabla \cdot \boldsymbol{J}_i, \quad \mbox{in} \ {\Omega _s},\ i = 1,2, \cdots ,K,\\
  \boldsymbol{J}_i &= -D_i\Big(\nabla c_i + \beta q_ic_i\nabla\phi  + \boldsymbol{N}_{k_i}(c_i)\Big), \notag
 \end{align}
where
\begin{align} \label{smpnp-nonliner}
 \boldsymbol{N}_{k_i}(c_i) =  k_{i}\frac{ c_{i} }{1-\sum\limits_{l} a_{l}^{3} c_{l}}\sum\limits_{l} a_{l}^{3} \nabla c_{l}, \ \
 k_i
 \left\{
   \begin{array}{l}
   = 0, \ {\rm for} \ {\rm NP\ equations}, \\
   \neq 0, \ {\rm for} \ {\rm SMNP\ equations}.
   \end{array} \right.
\end{align}

\subsection{Boundary conditions}

\begin{figure}[H]
\centerline{
 \includegraphics[scale=0.5]{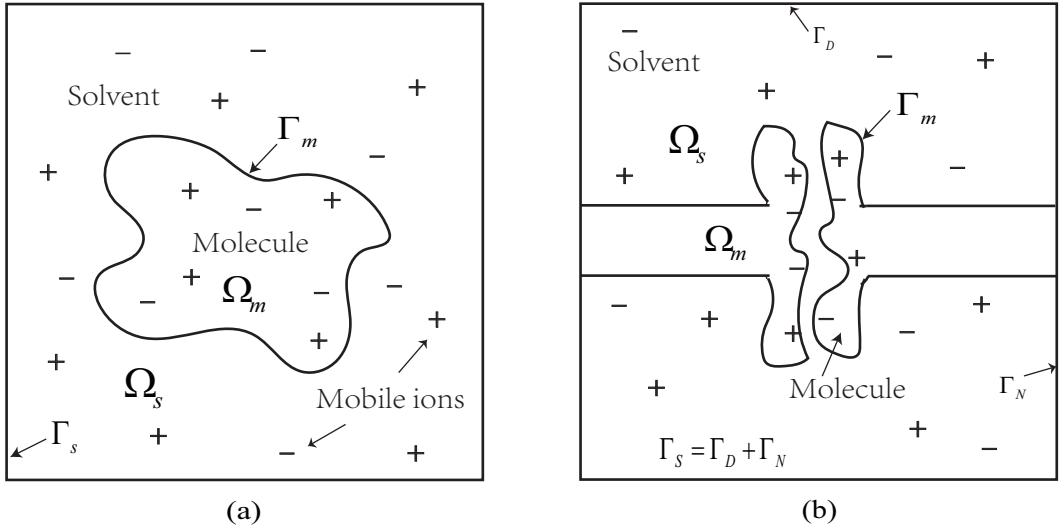} }
 \caption{ A 2-D schematic view of the biological systems: (a) a fixed biomolecule; (b) an ion channel (or similar a nanopore) embedded in a biomolecular membrane.}
 \label{region_smPNP_BD}
\end{figure}

\noindent
In this work, including the internal interface conditions \eqref{smpnp-poi-gamma} for the Poisson equation \eqref{smpnp-poi},  we consider the following boundary conditions
\begin{eqnarray} \label{smpnp-bd}
  \left\{\begin{array} {ll}
  [\phi]  = 0, \
  [\epsilon\frac{\partial\phi}{\partial\vec{n}} ]= 0, \   &\mbox{on} \ \Gamma_m,   \\
   \phi = \phi_0, \   &\mbox{on} \ \Gamma_s,      \\
   c_i = c_i^\infty, \ \ &\mbox{on} \ \Gamma_s     \\
   \boldsymbol{J}_i \cdot\vec{n} = 0, \  &\mbox{on} \ \Gamma_m,
  \end{array} \right.
 \end{eqnarray}
where $[\cdot]$ denotes the jump of the electrostatic potential at the internal interface $\Gamma_m$. The interface conditions \eqref{smpnp-poi-gamma} represent the continuity conditions for the electrostatic potential on the interface $\Gamma_m$. $\Gamma_s$  is the outer boundary of the solvent region $\Omega_s$,  in which the Dirichlet boundary $\Gamma_D$ and the Neumann boundary $\Gamma_N$ are all considered for the mixed boundary case. For example, domains and boundaries of demo systems are shown in  Fig.\ref{region_smPNP_BD}: a 2-D schematic view of biological systems. $c_i^\infty$ and $\phi_0$ are the bulk concentration of the $i$th ionic species and the applied potential, respectively. $\frac{\partial\phi}{\partial \vec{n} }$ denotes the normal derivative at the boundary with the exterior unit normal $\vec{n}$. The homogeneous Neumann boundary conditions preserve the conservation of the system and the continuity of the electrostatic potential at the internal interface $\Gamma_m$.

\subsection{Singular decomposition for the Poisson equation }

\noindent
In this paper, we only consider the steady-state PNP model, that is $\frac{\partial c_i}{\partial t} = 0$. To deal with the singular permanent charges, an effective strategy for solving Eq. \eqref{smpnp-poi} is to decompose the solution of the Poisson equation into three components: a singular component, a harmonic component and a regular component \cite{Chern2003Accurate,Lu2010pnp,Lu2011smpnp}, that is, $\phi=\phi^{s}+\phi^{h}+\phi^{r}$. For the sake of completeness, we introduce the decomposition process and their governing equations, respectively.

Firstly, the singular component $\phi^s$ is restricted into $\Omega_m$, and it is the solution of
\begin{eqnarray} \label{singular-part}
 - \epsilon_{m}\Delta \phi^{s}(x) = \rho^{f}(x), \quad x \in \mathbb{R}^{3}.
\end{eqnarray}
In fact, $\phi^{s}(x)$ can be given analytically by the sum of Coulomb potentials, that is
$$
\phi^{s}(x)=\sum_{j=1}^{N} \frac{q_{j}}{\epsilon_{m} 4 \pi\left|x-x_{j}\right|},
$$
where $N$ is the total number of particles in the biomolecule, and $|x-x_{j}|$ denotes the distance between the current position $x$ and the particle center $x_j$ of the $j$th ion species.

The harmonic component $\phi^{h}$ is the solution of a Laplace equation:
\begin{eqnarray} \label{harmonic-part}
\begin{array}{lr}
-\Delta \phi^{h}(x)=0, & x \in \Omega_{m}, \\
\phi^{h}(x)=-\phi^{s}(x), & x \in \Gamma_m .
\end{array}
\end{eqnarray}

Subtracting the above two components $\phi^{s}$ and $\phi^{h}$ from Eq.\eqref{smpnp-poi}, we get the governing equation of the regular component $\phi^{r}(x)$:
\begin{eqnarray} \label{regular-part}
-\nabla \cdot\left(\epsilon \nabla \phi^{r}(x, t)\right) = \lambda \sum_{i} q_{i} c_{i}(x, t), \quad x \in \Omega,
\end{eqnarray}
and the interface conditions
$$
\phi_{s}^{r}-\phi_{m}^{r}=0, \quad \epsilon_{s} \frac{\partial \phi_{s}^{r}}{\partial n}-\epsilon_{m} \frac{\partial \phi_{m}^{r}}{\partial n}=\epsilon_{m} \frac{\partial\left(\phi^{s}+\phi^{h}\right)}{\partial n}, \quad x \in \Gamma_m.
$$
It is worth noting that there is no decomposition of the electrostatic potential in the solvent region, thus $\phi(x)=\phi^{r}(x)$ in $\Omega_{s}$. Hence, the final regularized SMPNP/PNP equations consist of the regularized Poisson equation \eqref{regular-part} and the SMNP/NP equations
\begin{align} \label{sdt-smpnp-pnp}
  - \nabla\cdot \Big( D_i(x)(\nabla c_i(x) + \beta q_ic_i(x)\nabla\phi^r(x))  + \boldsymbol{N}_{k_i}(c_i(x)) \Big) = 0,  \quad x \in\ {\Omega _s}.
 \end{align}
 In the following content, we still also use $\phi$ to represent the regular component $\phi^{r}(x)$, and the singular and harmonic components have been considered to get the complete electrostatic potential inside molecules.

 Compared to the original model \eqref{smpnp-np}-\eqref{smpnp-poi-gamma}, the above decompositions (see \eqref{singular-part}-\eqref{regular-part}) have a number of nice properties. Firstly, the decomposition of the electrostatic potential only occurs inside biomolecules, so the numerical solution of $\phi^r$ in $\Omega_s$ does not possess the numerical instability problem \cite{Holst2012Adaptive}. Secondly, the singular and harmonic components only need to be solved one time in advance when decoupling the regularized SMPNP/PNP equations. More comments and relative comparisons can be fund in \cite{Lu2010pnp} and the references therein.

\subsection{A transformed form of the SMPNP equations }

\noindent
By introducing a set of generalized Slotboom variables \cite{Lu2011smpnp}, the regularized SMPNP equations \eqref{regular-part}-\eqref{sdt-smpnp-pnp} can be written as
\begin{align}
 - \nabla \cdot\left(\bar{D_i} \nabla \bar{C_i}\right) = 0, \label{smpnp-slotboom-np}\\
 -\nabla\cdot\big(\epsilon\nabla\phi\big) - \lambda\sum\limits_i^K q_i\bar{C_i} e^{-\Psi_{i}} = 0, \label{smpnp-slotboom-poi}
\end{align}
with
\begin{equation} \label{smpnp-slotboom}
\left\{\begin{array} {rcl}
\Psi_{i}  &=& \beta q_{i} \phi - k_{i} \ln \Big(1-\sum\limits_{l}^{K} a_{l}^{3} c_{l}\Big), \ \
              k_{i}=a_{i}^{3} / a_{0}^{3}, \\
\bar{D_i} &=& D_{i} e^{-\Psi_{i}}, \\
\bar{C_i} &=& c_{i} e^{\Psi_{i}}.
\end{array}  \right.
\end{equation}
 Physically, $\Psi_{i}$ can be seen as a modification of the electrostatic potential $\phi$ due to the size effects. If the size effect is not considered ($k_i = 0$ or $a_i = 0$), the transformed forms \eqref{smpnp-slotboom-np}-\eqref{smpnp-slotboom} reduce to the classical Slotboom transform of PNP equations \cite{Lu2010pnp,Slotboom1973Computer,Tu2013A}.

The transformed SMNP equation \eqref{smpnp-slotboom-np} is a self-adjoint second-order elliptic equation about the Slotboom variable $\bar{C_i}$. Different from the classical NP equations, the coefficient $\bar{D_i}$ depends on $\phi$ and $c_{i}$ in SMNP equations. Therefore, a semi-implicit scheme is employed in our scheme. In the iterative process of equations decoupling, we use the solution at the $(n-1)$th step $c_i^{n-1}$ to calculate the coefficient $\bar{D_i}$, and then solve the transformed Eq. \eqref{smpnp-slotboom-np} to obtain the solution at the current $n$th step $c_i^{n}$. This strategy can make the stiffness matrices symmetric for the generalized Slotboom variable $\bar{C_i}$. And the condition number of the stiffness matrix derived from discretizing the transformed Eq. \eqref{smpnp-slotboom-np} may be smaller than that of the stiffness matrix produced by the origianl Eq. \eqref{sdt-smpnp-pnp}. Thus the decoupling iterative methods applied to the linear system might converge faster \cite{Lu2011smpnp}. However, in biomolecular simulations, as shown in \cite{Lu2010pnp}, the discretization of the transformed Eq. \eqref{smpnp-slotboom-np}  always leads to an ill-conditioned stiffness matrix because a strong electrostatic field exists near the molecular surface. In addition, the introduction of Slotboom variables makes the Poisson equation \eqref{smpnp-slotboom-poi} become nonlinear for the electrostatic potential $\phi$. So a nonlinear iterative scheme, e.g. Newton method, is necessary for solving the nonlinear Poisson equation, which may cost much more CPU time. In practical numerical simulations, our previous experience shows that the Newton method is sensitive to the initial value, especially in macromolecular biological channel simulations, e.g. KcsA ${\rm K^+}$ channels, see \cite{Liu2017Incorporating,Wang2021A}. In order to avoid multi-level nonlinear iteration and improve the efficiency and robustness of the our method, we use the normal unknown variables when solving SMPNP equations \eqref{smpnp-slotboom-np}-\eqref{smpnp-slotboom-poi} in this work.

At first, we use the the normal variables $c_i$ to rewrite  \eqref{smpnp-slotboom-np}-\eqref{smpnp-slotboom} as follows:
\begin{align}
 - \nabla \cdot\left(D_{i} e^{-\Psi_{i}}\nabla (e^{\Psi_{i}}c_{i}) \right) = 0, & \ \ {\rm in} \ \Omega_s,
 \label{smpnp-slotboom-np-norm}\\
 -\nabla\cdot\big(\epsilon\nabla\phi\big) = \lambda\sum\limits_i^K q_ic_i, & \ \ {\rm in} \ \Omega, \label{smpnp-slotboom-poi-norm}
\end{align}
with
\begin{equation} \label{smpnp-slotboom-Psi}
\Psi_{i} = \beta q_{i} \phi - k_{i} \ln \Big(1-\sum\limits_{l}^{K} a_{l}^{3} c_{l}\Big), \ \
k_{i} = a_{i}^{3} / a_{0}^{3}.
\end{equation}
Let $u = \beta e_c\phi$ to nondimensionalize the electrostatic potential, and $q_i = z_i e_c$. Then Eqs. \eqref{smpnp-slotboom-np-norm}-\eqref{smpnp-slotboom-Psi} become:
\begin{align}
 - \nabla \cdot\left(D_{i} e^{-\Psi_{i}}\nabla (e^{\Psi_{i}}c_{i}) \right) = 0, & \ \ {\rm in} \ \Omega_s,
 \label{smpnp-slotboom-np-nonnorm}\\
 -\nabla\cdot\big(\epsilon\nabla u\big)
 =\beta e_{c}^{2} \lambda\sum\limits_i^K z_ic_i, & \ \ {\rm in} \ \Omega, \label{smpnp-slotboom-poi-nonnorm}
\end{align}
where
\begin{equation} \label{smpnp-slotboom-nonPsi}
 \Psi_{i} = z_iu - k_i  \ln \left(1-\sum\limits_{l}^{K} a_{l}^{3} c_{l}\right), \ \
 k_{i} = a_{i}^{3} / a_{0}^{3}.
\end{equation}
The corresponding boundary conditions become
\begin{eqnarray} \label{sdt-smpnp-nonbd-exp-c}
  \left\{\begin{array} {ll}
  [u] = 0, \ [\epsilon\frac{\partial u}{\partial\vec{n} }] &= \beta e_c \epsilon_{m} \frac{\partial\left(\phi^{s}+\phi^{h}\right)}{\partial \vec{n}}, \ \mbox{on} \ \Gamma_m,  \\
   u = u_0,  \ &\mbox{on} \ \Gamma_D,     \\
   c_i = c_i^\infty, \ &\mbox{on} \ \Gamma_D    \\
   \boldsymbol{J}_i\cdot\vec{n} = 0, \  &\mbox{on} \ \Gamma_m,
  \end{array} \right.
 \end{eqnarray}
where
$ u_0 = \beta e_c\phi_0 $,
$\boldsymbol{J}_i = -\left(D_{i} e^{-\Psi_i} \nabla (e^{\Psi_i}c_{i}) \right)$, and $\phi_0$, $c_{i}^\infty$ are the applied potential and bulk concentrations defined by \eqref{smpnp-bd}.

\section{Inverse averaging finite element method for SMPNP equations }
\label{fem-eafe}

\noindent
In this section, we will introduce the IAFEM for the reformulated SMPNP equations \eqref{smpnp-slotboom-np-nonnorm}-\eqref{sdt-smpnp-nonbd-exp-c} in detail. In order to facilitate the presentation and understanding of the latter contents, we first report some preliminary notations of the finite element discretization.  Let $H^1(\Omega)$ be the Sobolev space of weakly differentiable functions. Denoted by
$$
\begin{aligned}
L^2(\Omega) & \equiv\left\{w: \Omega \rightarrow \mathbb{R} \mid \int_{\Omega} w^{2} d x<\infty\right\}, \\
H^1(\Omega) & \equiv\left\{w \in L^2(\Omega)\mid \int_{\Omega}|\nabla w|^{2} d x<\infty\right\}, \\
H_0^1(\Omega) & \equiv\left\{w \in H^1(\Omega)\mid w|_{\Gamma_{D}}=0\right\}
\end{aligned}
$$
be the spaces with associated norms:
$$
\|w\|_{0} \equiv\left(\int_{\Omega} w^{2} d x\right)^{1 / 2}, \quad|w|_{1} \equiv\left(\int_{\Omega}|\nabla w|^{2} d x\right)^{1 / 2}, \quad\|w\|_{1}^{2} \equiv|w|_{1}^{2}+\|w\|_{0}^{2},
$$
and $\|w\|_{\infty} \equiv \sup _{x \in \Omega}|w(x)|$, and the inner product $(f,g)_{\Omega} := \int_\Omega f g dx$.

\subsection{Weak forms and finite element discretization of the reformulated SMPNP equations}
\noindent
Integrating by parts, and noting the interface conditions in \eqref{sdt-smpnp-nonbd-exp-c}, the weak forms of the reformulations \eqref{smpnp-slotboom-np-nonnorm}-\eqref{sdt-smpnp-nonbd-exp-c} are to find $c_i \in H^1(\Omega_s)$ ($1\le i \le K$) and $u \in H^1(\Omega)$ satisfying
\begin{align} \label{std-smpnp-NP-weak}
 \big( D_{i} e^{-\Psi_i} \nabla (e^{\Psi_i}c_{i}), \nabla v \big)_{\Omega_s} = 0, \ \ \forall v \in H_0^1(\Omega), \\
 \big( \epsilon\nabla u, \nabla w \big)_{\Omega} = \Big(\beta e_{c}^{2} \lambda\sum\limits_i^K z_ic_i,  w \Big)_{\Omega} - \beta e_c\Big(\epsilon_{m} \frac{\partial\left(\phi^{s}+\phi^{h}\right)}{\partial \vec{n}}, w\Big)_{\Gamma_m}, \ \ \forall w \in H_0^1(\Omega), \label{std-smpnp-Poi-weak}
\end{align}
where $\Psi_i$ is defined by \eqref{smpnp-slotboom-nonPsi}.

Let $\mathcal{T}_{h} = \{ T \}$ be a triangulation of $\Omega$ with (triangular/tetrahedral) elements $T$,  $X_{h}=\left\{q_{i}\right\}_{i=1}^{N_{v}}$ be the set of all vertices of
$\mathcal{T}_{h}$.
Let $V_h \subset H_0^1(\Omega)$ be the piecewise linear finite element space, and $H_0^1(\Omega)$ is a Sobolev space of weakly differentiable functions which vanish on the boundary of the domain $\Omega$. Denote the nodal basis function in $V_{h}$ with $\varphi_{i},\ i = 1, 2, \cdots, N_v$, which is linear on the $T$ and
\begin{align} \label{Lagrange-basis}
 \varphi_{i}\left(q_{i}\right)=1, \quad \varphi_{i}\left(q_{j}\right)=0, \quad j \neq i.
\end{align}
For a given $T \in \mathcal{T}_{h}$, we have
\begin{align} \label{element-stiff-identity-0}
 \int_{T} \nabla u_{h} \cdot \nabla v_{h} d x=\sum_{i, j} e_{i j}^{T} u_{h}\left(q_{i}\right) v_{h}\left(q_{j}\right), \quad \forall u_{h}, v_{h} \in V_{h}.
\end{align}
Note that $e_{i j}^{T}=\int_{T} \nabla \varphi_{j} \cdot \nabla \varphi_{i} d T$ represents some geometric information of the element $T,$ and it holds for linear Lagrangian finite element basis functions that $e_{i i}^{T}=-\sum\limits_{j \neq i} e_{i j}^{T}$. Then we can easily transform \eqref{element-stiff-identity-0} to the following simple but important identity
\begin{align} \label{element-stiff-identity-1}
 \int_{T} \nabla u_{h} \cdot \nabla v_{h} dx =
 -\sum_{i,j} e_{i j}^{T}(u_{h}(q_{i})-u_{h}(q_{j}))(v_{h}(q_{i})-v_{h}(q_{j})), \  \forall u_{h}, v_{h} \in V_{h}.
\end{align}
Let the test function $v_h$ take the linear Lagrangian basis function $\varphi_i$ at $q_i$, we have
\begin{align} \label{element-stiff-identity-2}
 \int_{T} \nabla u_{h} \cdot \nabla \varphi_i dx =
 -\sum_{q_j\in T} e_{i j}^{T}(u_{h}(q_{i})-u_{h}(q_{j})).
\end{align}

The finite element discretization of \eqref{std-smpnp-NP-weak}-\eqref{std-smpnp-Poi-weak} is to find $c_{i,h} \in V_h(\Omega_s),\ 1\le i \le K$, and $u_h \in V_h(\Omega)$, such that
 \begin{align} \label{std-smpnp-NP-FEMweak}
  \big( D_{i} e^{-\Psi_{i,h}} \nabla (e^{\Psi_{i,h}}c_{i,h}), \nabla v_h \big)_{\Omega_s} = 0, \ \ \forall v_h \in V_h, \\
  \big( \epsilon\nabla u_h, \nabla w_h \big)_{\Omega} = \Big(\beta e_{c}^{2} \lambda\sum\limits_{i=1}^K z_ic_{i,h},  w_h \Big)_{\Omega} - \beta e_c\Big(\epsilon_{m} \frac{\partial\left(\phi^{s}+\phi^{h}\right)}{\partial \vec{n}}, w_h \Big)_{\Gamma_m}, \ \ \forall w_h \in V_h, \label{std-smpnp-Poi-FEMweak}
 \end{align}
 where
 \begin{equation} \label{smpnp-slotboom-FEMnonPsi}
  \Psi_{i,h} = z_iu_h - k_i  \ln \left(1-\sum\limits_{l}^{K} a_{l}^{3} c_{l,h}\right), \ \
  k_{i} = a_{i}^{3} / a_{0}^{3}.
 \end{equation}

In order to show differences between finite element approximations of the reformulated and traditional schemes clearly, we also present the standard finite element scheme for \eqref{sdt-smpnp-pnp} as follows:

For each $i$, $1\le i \le K$, find $c_{i,h} \in V_h(\Omega_s)$ and $u_h \in v_h(\Omega)$, such that
 \begin{align}
  \left(D_{i}\Big(\nabla c_{i,h} + z_ic_{i,h}\nabla u_h + \frac{k_{i} c_{i,h} }{1-\sum\limits_{l} a_{l}^{3} c_{l,h} }\sum\limits_{l} a_{l}^{3} \nabla c_{l,h}  \Big), \nabla v_h\right)_{\Omega_s} = 0, \ \ \forall v_h \in V_h, \label{std-smpnp-np-stdFE}
 \end{align}
and the finite element discredization for Poisson equation is same as \eqref{std-smpnp-Poi-FEMweak}.\\

We decouple the nonlinear coupling system \eqref{std-smpnp-NP-FEMweak}-\eqref{std-smpnp-Poi-FEMweak} with Gummel iteration \cite{Gummel1964A}. In each iteration, the Poisson equation and each NP equation are solved successively. The ion concentrations are treated as known functions when solving the electrostatic potential, and vice versa. The process repeats until the difference of solutions in two adjacent iterations becomes smaller than a given tolerance.

For the steady-state case, in order to make iterations between the Poisson and SMNP equations converge, it is necessary to employ the under-relaxation technique, especially when macromolecules exist. In other words, solutions are updated with a linear combination of solutions respectively obtained from the last iteration and the current iteration, rather than just using solutions derived from the current iteration. This under-relaxation scheme \cite{Lu2011smpnp,Lu2007Electrodiffusion,Tu2013A} is described by
\begin{align*}
 &u^{\text {new }}=\alpha u^{\text {old }}+(1-\alpha) u^{\text {new }}, \\
 &c_{i}^{\text {new }}=\alpha c_{i}^{\text {old }}+(1-\alpha) c_{i}^{\text {new }}, \quad i=1, 2, \cdots, K,
\end{align*}
where the relaxation parameter $0 < \alpha <1$ is a predefined constant. We note that without the under-relaxation technique, the iterations may not converge.  More specifically, the iterative process will be presented in Section \ref{numerical-experiment}.

\subsection{Inverse averaging technique}

\noindent
From the aforementioned content, if one wants to solve the finite element approximation  equation \eqref{std-smpnp-NP-FEMweak} with the normal unknown variable $c_i$ accurately, the numerical difficulty lies in dealing with the exponential coefficients $e^{-\Psi_i}$ and $e^{\Psi_i}$. In this subsection, we apply a novel inverse averaging technique which calculates the inverse averages of the exponential coefficients on the edge $E_{i j}=\overline{q_{i} q_{j}}$ of the element $T$.

\subsubsection{The inverse average of the exponential coefficient }

\noindent
The inverse average of the exponential coefficient on the edge $E_{i j}=\overline{q_{i} q_{j}}$ is denoted with
\begin{align} \label{inverse-average-exp}
 E(\Psi)_{E_{i j}}=\left(\frac{\int_{q_{i}}^{q_{j}} e^{- \Psi} \mathrm{d} s}{\left|E_{i j}\right|}\right)^{- 1} \triangleq I\left(E_{i j}\right),
\end{align}
where $\Psi$ is defined by \eqref{smpnp-slotboom-nonPsi} in this paper. In the other cases, $\Psi$ may have its own specific definition, e.g. \cite{ZhangQianru2021A,Zhanng2021An}.
The work \cite{ZhangQianru2021A} shows that the inverse averaging technique is significant to control the effect of large electrostatic fields on currents and enhance the stability of numerical methods for solving the standard PNP equations, particularly with rapidly varying coefficients when solving the three-dimensional drift-diffusion model in semiconductor device simulations. 
Inspired by this, we also use the inverse averaging technique to deal with the exponential coefficients in our work for solving the reformulated finite element approximation Eq. \eqref{std-smpnp-NP-FEMweak} when simulating biomolecular systems and ion channels.

Firstly, similar to \cite{ZhangQianru2021A}, we also assume that $\Psi$ is linear on the edge $E_{i j}$, that is
 \begin{align} \label{Psi-newA4}
  \Psi(\mathbf{x})=\left(\frac{\Psi_{j}-\Psi_{i}}{\left|E_{i j}\right|}\right)\left(\mathbf{x}-\mathbf{x}_{q_{i}}\right)+\Psi_{i}, \quad \mathbf{x} \in\left[\mathbf{x}_{q_{i}}, \mathbf{x}_{q_{j}}\right].
 \end{align}
From \eqref{inverse-average-exp} and \eqref{Psi-newA4}, we get
\begin{align}\label{IE-newA4}
  I(E_{i j})
  &= \left(\int_{q_{i}}^{q_{j}} \frac{e^{-\Psi_{i}}}{\left|E_{ij}\right|}
    \left(\frac{e^{\Psi_{i}}}{e^{\Psi_{j}}}\right)^{\frac{\mathbf{x}-\mathbf{x}_{q_{i}}}{|E_{i j}|}} d\mathbf{x}\right)^{-1}
   = \left(\int_{q_{i}}^{q_{j}} \frac{e^{-\Psi_{i}}}{\left|E_{ij}\right|}
    e^{(\Psi_{i}-\Psi_{j}) \cdot {\frac{\mathbf{x}-\mathbf{x}_{q_{i}}}{|E_{i j}|}} } d\mathbf{x} \right)^{-1} \notag\\
  &= \left( \frac{e^{-\Psi_{i}}}{\Psi_{i}-\Psi_{j}} \int_{q_{i}}^{q_{j}}
    e^{(\Psi_{i}-\Psi_{j}) \cdot {\frac{\mathbf{x}-\mathbf{x}_{q_{i}}}{|E_{i j}|}} } d\Big( {\frac{(\Psi_{i}-\Psi_{j}) (\mathbf{x}-\mathbf{x}_{q_{i}})}{|E_{i j}|}} \Big) \right)^{-1} \notag\\
  &= \left( \frac{e^{-\Psi_{i}}}{\Psi_{i}-\Psi_{j}} e^{(\Psi_{i}-\Psi_{j}) \cdot {\frac{\mathbf{x}-\mathbf{x}_{q_{i}}}{|E_{i j}|}} } \Big|_{x_{q_i}}^{x_{q_j}}  \right)^{-1} \notag\\
  &=e^{\Psi_{i}} B\left(\Psi_{i}-\Psi_{j}\right),
\end{align}
 where $B(t)$ is the Bernoulli function defined by
 $$
 B(t)=\left\{\begin{array}{ll}
 \frac{t}{e^{t}-1}, & t \neq 0, \\
 1, & t=0 .
 \end{array}\right.
 $$
For numerical stability, if the difference between two nodal values of $\Psi_{i}$ is very small, the corresponding terms ``$ B(\Psi_i - \Psi_j)$'' should be calculated using Taylor expansions. For more details, please refer to \cite{ZhangQianru2021A} and the references therein. Especially, the Bernoulli function $B(t)$ is calculated by
$$
 B(t)=\left\{\begin{array}{ll}
 \frac{t}{e^{t}-1}, & |t| > 10^{-4}, \\
  \Big( \big(-\frac{1}{720} t^2 + \frac{1}{12} \big) t - \frac{1}{2} \Big) t + 1, & \text{otherwise},
 \end{array}\right.
$$
in our computation.

\subsubsection{The inverse averaging finite element scheme for the reformulated SMNP equations}

\noindent
 Now we give a derivation of the IAFEM for the reformulated Eq. \eqref{std-smpnp-NP-FEMweak}. First of all, referring to the mean value theorem of integrals and approximating the exponential coefficient $e^{-\Psi}$ with $E(-\Psi)_{E_{ij}}$  on the edge $E_{ij}$ of the element $T$, we have
\begin{align} \label{std-smpnp-NP-FEMweak-int}
 0 &= \big( D_{i} e^{-\Psi_{i,h}} \nabla (e^{\Psi_{i,h}}c_{i,h}), \nabla v_h \big)_{\Omega_s} \notag\\
   &= \int_{\Omega_s} D_{i} e^{-\Psi_{i,h}} \nabla (e^{\Psi_{i,h}}c_{i,h}) \cdot \nabla v_h d\Omega_s \notag\\
   &= \sum\limits_{T\in T_h} \int_T D_{i} e^{-\Psi_{i,h}} \nabla (e^{\Psi_{i,h}}c_{i,h}) \cdot \nabla v_h dT \notag\\
   &\approx \sum\limits_{T\in T_h} D_{i} E(-\Psi)_{E_{ij}} \int_T \nabla (e^{\Psi_{i,h}}c_{i,h}) \cdot \nabla v_h dT.
\end{align}
In the following content, we introduce the computation of the element-wise stiffness matrix for \eqref{std-smpnp-NP-FEMweak-int}, i.e., $A = (a_{i j}^{T})_{T \in \mathcal{T}_{h}}$, in detail.

On a element $T$, let $v_{h}$ take the associated piecewise linear finite element basis function. By using \eqref{element-stiff-identity-2},  we have
\begin{align} \label{element-stiffness-0}
 &D_{i} E(-\Psi)_{E_{ij}} \int_T \nabla (e^{\Psi_{i,h}}c_{i,h}) \cdot \nabla v_h dT  \notag\\
 &= D_{i} E(-\Psi)_{E_{ij}} \sum\limits_{q_j\in T}(e^{\Psi_{i,h}}c_{i,h})(q_j) \int_T \nabla\varphi_j \cdot \nabla\varphi_i dT  \notag\\
 &= D_{i} E(-\Psi)_{E_{ij}} \sum\limits_{q_j\in T}(e^{\Psi_{i,h}}c_{i,h})(q_j) e_{ij}^T  \notag\\
 &= - D_{i}\sum\limits_{q_j\in T, q_j\ne q_i} E(-\Psi)_{E_{ij}}
    \big((e^{\Psi_{i,h}}c_{i,h})(q_i) - (e^{\Psi_{i,h}}c_{i,h})(q_j) \big) e_{ij}^T.
\end{align}
Furthermore, the approximated coefficient $E(-\Psi)_{E_{ij}}$ is calculated by \eqref{IE-newA4} on the edge $E_{ij}$, then we get
\begin{align} \label{element-stiffness-1}
 &D_{i} E(-\Psi)_{E_{ij}} \int_T \nabla (e^{\Psi_{i,h}}c_{i,h}) \cdot \nabla v_h dT  \notag\\
 &= - D_{i}\sum\limits_{q_j\in T, q_j\ne q_i} E(-\Psi)_{E_{ij}}
    \big((e^{\Psi_{i,h}}c_{i,h})(q_i) - (e^{\Psi_{i,h}}c_{i,h})(q_j) \big) e_{ij}^T  \notag\\
 &= - \sum\limits_{q_j\in T, q_j\ne q_i} D_i e^{-\Psi_{i,h}} B(\Psi_j-\Psi_i)(e^{\Psi_{i,h}}c_{i,h})(q_i) e_{ij}^T
   + \sum\limits_{q_j\in T, q_j\ne q_i} D_i (e^{-\Psi_{i,h}} B(\Psi_i-\Psi_j)(e^{\Psi_{i,h}}c_{i,h})(q_j) e_{ij}^T \notag\\
 &= \left( - \sum\limits_{q_j\in T, q_j\ne q_i} D_i B(\Psi_j-\Psi_i) e_{ij}^T \right) c_{i,h}(q_i)
     + \sum\limits_{q_j\in T, q_j\ne q_i} \Big( D_i B(\Psi_i-\Psi_j)e_{ij}^T \Big) c_{i,h}(q_j).
\end{align}
The nonzero entries of the element-wise stiffness matrix $A = (a_{i j}^{T})_{T \in \mathcal{T}_{h}}$ can be written as
$$
a_{i j}^{T}=
\begin{cases}
 D_{i} B(\Psi_{i}-\Psi_{j}) e_{i j}^{T}, & j \neq i, \\
 -\sum\limits_{k \neq i} D_{i} B(\Psi_{k}-\Psi_{i}) e_{i k}^{T}, & j=i,
\end{cases}
$$
where
  $$
  e_{i j}^{T}=\int_{T} \nabla \varphi_{j} \cdot \nabla \varphi_{i} d T,
  $$
and $\varphi_{i}, i=1, \ldots N_{h}$ are nodal basis functions in $V_{h}$, which satisfy \eqref{Lagrange-basis}. \\

$~$\\

\begin{algorithm} [H]
\caption{Gummel iteration for SMPNP equations with IAFEM} \label{algo-gummel-IEAFE}

 Step 1: Initialization for nonlinear iteration:

 \quad\quad Initialize error tolerance $tol$ and maximize iteration number $max_{-}N$ ;

 \quad\quad Initialize iterative step counter $n = 0$;

 \quad\quad Initialize electrostatic potential $u_h^n = 0$;

 \quad\quad for each $i \in 1, 2, \cdots, K$ do

 \quad\quad\quad\quad Initialize concentration $c_{i,h}^{n}=0$ and auxiliary variable $c_{i,h}^{n+1}=0$;

 \quad\quad end for \\

 Step 2: Nonlinear iteration: $n \ge 1$, solving the decoupled equations:

 \While {$\left\|u_h^{n+1}-u_h^{n}\right\|> tol$ and $j \le max_{-}N$}
 {
    \begin{align}
     &\big( \epsilon\nabla u_h^{n+1}, \nabla w_h \big)_{\Omega} = \Big(\beta e_{c}^{2} \lambda\sum\limits_i^K z_ic_{i,h}^{n+1},  w_h \Big)_{\Omega} - \beta e_c\Big(\epsilon_{m} \frac{\partial\left(\phi^{s}+\phi^{h}\right)}{\partial \vec{n}}, w_h \Big)_{\Gamma_m},  \label{algo-gummel-poi}\\
     &\Psi_{i,h}^{n} = z_iu_h^{n} - k_i  \ln \left(1-\sum\limits_{l}^{K} a_{l}^{3} c_{l,h}^{n}\right), \ \  k_{i} = a_{i}^{3} / a_{0}^{3},  \label{algo-gummel-Psi}\\
     &\big( D_{i} e^{-\Psi_{i,h}^{n}} \nabla (e^{\Psi_{i,h}^{n}}c_{i,h}^{n+1}), \nabla v_h \big)_{\Omega_s} = 0,  \label{algo-gummel-np}
    \end{align}

    \eIf {$\|u_h^{n+1}-u_h^{n}\|< tol$}
     {
       break;
     }
     {
       $u_h^{n+1} \leftarrow \alpha u_h^{n}+(1-\alpha) u_h^{n+1} $;

      $c_{i,h}^{n+1} \leftarrow \alpha c_{i,h}^{n}+(1-\alpha) c_{i,h}^{n+1} $;
     }
     $ n \leftarrow n + 1$;
 }

 Step 3: Output the electrostatic potential $\phi_h^{n+1} :=u_h^{n+1} / (\beta e_c)$ and concentrations $c_{i,h}^{n+1}$.
\end{algorithm}

\setcounter{equation}{0}
\section{Numerical tests and application in simulating biomolecular systems and ion channels} \label{numerical-experiment}

\noindent
In this section, we will use the IAFEM to solve the SMPNP (PNP) equations for simulating biomolecular systems and ion channels. To demonstrate the accuracy and robustness of the IAFEM, a model problem with analytic solutions on a cube is firstly tested. Then the size-effect simulations on a biomolecular sphere with different charges and an ion channel are respectively implemented. All the numerical algorithms are implemented based on the three-dimensional parallel finite element toolbox Parallel Hierarchical Grid (PHG) \cite{Zhang2009A}. The computations were done on the high performance computers of State Key Laboratory of Scientific and Engineering Computing, Chinese Academy of Sciences.\\

As mentioned above, the Gummel iterative method is used to decouple the nonlinear coupling system \eqref{std-smpnp-NP-FEMweak}-\eqref{std-smpnp-Poi-FEMweak}, and also used for solving the traditional finite element approximation \eqref{std-smpnp-Poi-FEMweak} and \eqref{std-smpnp-np-stdFE}. In order to clearly understand the Gummel iterative process and calculation process in this work, we present the iterative process in Algorithm \ref{algo-gummel-IEAFE}.

Similarly, the nonlinear iterative process of the traditional FEM for SMPNP equations is presented as follows: \\

\begin{algorithm} [H]
 \caption{Gummel iteration for SMPNP equations by FE}\label{algo-gummel-FE}

 Step 1: The same as step 1 in Algorithm \ref{algo-gummel-IEAFE};

 Step 2: Nonlinear iteration: $n\ge 1$, solving the coupled equations:

\While {$\left\|u_h^{n+1}-u_h^{n}\right\|> tol$ and $n \le max_{-}N$}
 {
  Step 2.1 Solving Poisson equation for $u_h^{n+1}$ through Eq. \eqref{algo-gummel-poi};

  Step 2.2 Solving SMNP equations for $c_{i,h}^{n+1}$ through the following equation
  \begin{align}
   &\left(D_{i}\Big(\nabla c_{i,h}^{n+1} + z_ic_{i,h}^{n+1}\nabla u_h^n + \frac{k_{i} c_{i,h}^{n+1} }{1-\sum\limits_{l} a_{l}^{3} c_{l,h}^{n} }\sum\limits_{l} a_{l}^{3} \nabla c_{l,h}^{n}  \Big), \nabla v_h\right)_{\Omega_s} = 0, \label{algo-gummel-np-FE}
  \end{align}
  \eIf {$\left\|u^{n+1}-u^{n}\right\|< tol$}
     {
       break;
     }
     {
       $u_h^{n+1} \leftarrow \alpha u_h^{n}+(1-\alpha) u_h^{n+1} $;

      $c_{i,h}^{n+1} \leftarrow \alpha c_{i,h}^{n}+(1-\alpha) c_{i,h}^{n+1} $;
     }
     $ n \leftarrow n + 1$;
 }

 Step 3: Output the electrostatic potential $\phi_h^{n+1} :=u_h^{n+1} / (\beta e_c)$ and concentrations $c_{i,h}^{n+1}$.
\end{algorithm}

\clearpage

 In the following context, some numerical experiments are reported to verify the effectiveness and robustness of the IAEEM for solving SMPNP equations. At first, a model problem with analytic solutions is presented to numerically verify the accuracy and order of the new scheme. Then, a sphere model and an ion channel system are separately considered.

\subsection{Accuracy and convergence tests }

\noindent
In this subsection, we report the numerical accuracy tests on the IAFEM for SMPNP equations through a model problem with analytic solutions.

\begin{example} \label{smpnp-test-sin}
 In this example, we consider a model problem with analytic solutions on a cube. Let the computational domain $\Omega =[0\AA, 1\AA]^3 $, and two charged species ${\rm K^+}$ and ${\rm Cl^-}$ are considered in the system. Specially, we use $c_p$ and $c_n$ to denote ${\rm K^+}$ and ${\rm Cl^-}$  concentrations only in this test, and their diffusion coefficients are $D_p = D_{\rm K^+} = 0.196 ~\AA^2/ps$, $D_n = D_{\rm Cl^-} = 0.203 ~\AA^2/ps$.
 Consider the following dimensionless SMPNP model problem
 \begin{eqnarray} \label{example-smpnp-1}
   \begin{cases}
    -\nabla \cdot(\nabla u) = (c_p-c_n) + f_u, & {\rm in} \ \Omega,\\
    -\nabla \cdot D_p(\nabla c_p + c_p \nabla u + \frac{k_{p} c_{p} }{1- \gamma\sum\limits_{l} a_{l}^{3} c_{l}}\sum\limits_{l} a_{l}^{3} \nabla c_{l}) = f_p, & {\rm in} \ \Omega,\\
    - \nabla \cdot D_n(\nabla c_n - c_n \nabla u + \frac{k_{n} c_{n} }{1- \gamma\sum\limits_{l} a_{l}^{3} c_{l}}\sum\limits_{l} a_{l}^{3} \nabla c_{l} ) = f_n, &{\rm in} \ \Omega,
   \end{cases}
 \end{eqnarray}
 where $k_{i}=a_{i}^{3} / a_{0}^{3}$ \ $(i=p,~n)$, the solvent molecular size $a_0 = 3.1 \AA$, the ion size $a_p = 1.51 \AA$, $a_n = 2.37 \AA$, and $\gamma = 6.022140857 \times 10^{-4} $. The right-hand functions $f_u$,  $f_p$, $f_n$ and the boundary conditions are respectively given by the following analytic solutions
 \begin{eqnarray} \label{exact-m3d-sin}
  \begin{cases}
   u = \sin(\pi x)\sin(\pi y)\sin(\pi z), \\
   c_{p} = \sin(2\pi x)\sin(2\pi y)\sin(2\pi z), \\
   c_{n} = \sin(3\pi x)\sin(3\pi y)\sin(3 \pi z).
  \end{cases}
 \end{eqnarray}
\end{example}

In this example, the piecewise linear finite element basis functions are used to
discretize the SMPNP model \eqref{example-smpnp-1}. The $L^2$ norm and $H^1$ norm errors are reported in Table \ref{3dsmpnp-L2error-H1error-order-sin}. The first column represents the mesh size of the uniform meshes. The numerical errors in  $L^2$ norm and $H^1$ norm are second-order and first-order reduction, respectively. This numerically demonstrates the convergence accuracy and reliability of the IAFEM for solving SMPNP equations.

\begin{table}[H]
\centering
\caption{ $L^2$ and  $H^1$ errors of the IAFEM for SMPNPEs. (Example \ref{smpnp-test-sin}) }
\vskip 0.1cm
\begin{tabular}{p{1.5cm}p{2.2cm}p{1.5cm}p{2.2cm}p{1.5cm}p{2.2cm}p{.9cm}}
 \hline\noalign{\smallskip}
 $h$ &$\|u_h-u\|_0$ &Order& $\|c_{p,h}-c_p\|_0$ & Order & $\|c_{n,h}-c_n\|_0$ & Order  \\
 \noalign{\smallskip}\midrule[1pt]\noalign{\smallskip}
  1/4   & 6.1793E-02 & $-$  & 8.0782E-02 & $-$  & 1.3742E-01 & $-$    \\
  1/8   & 1.9120E-02 & 1.69 & 6.2609E-02 & 0.37 & 9.6801E-02 & 0.51   \\
  1/16  & 5.0597E-03 & 1.92 & 1.9858E-02 & 1.66 & 3.8814E-02 & 1.32   \\
  1/32  & 1.2967E-03 & 1.96 & 5.2659E-03 & 1.91 & 1.1011E-02 & 1.82   \\
  1/64  & 3.3956E-04 & 1.93 & 1.3353E-03 & 1.98 & 2.9460E-03 & 1.90   \\
  1/128 & 8.7309E-05 & 1.96 & 3.3437E-04 & 2.00 & 7.4430E-03 & 1.98   \\
 \noalign{\smallskip}\hline
 \hline\noalign{\smallskip}    
 $h$ &$\|u_h-u\|_1$ &Order& $\|c_{p,h}-c_p\|_1$ & Order & $\|c_{n,h}-c_n\|_1$ & Order  \\
 \noalign{\smallskip}\midrule[1pt]\noalign{\smallskip}
  1/4   & 4.5797E-01 & $-$  & 1.0250E+00 & $-$  & 3.1852E+00 & $-$    \\
  1/8   & 1.8227E-01 & 1.33 & 9.0683E-01 & 0.18 & 2.2351E+00 & 0.51   \\
  1/16  & 8.1152E-02 & 1.17 & 3.5981E-01 & 1.33 & 9.2410E-01 & 1.27   \\
  1/32  & 3.9094E-02 & 1.05 & 1.5825E-01 & 1.19 & 3.8890E-01 & 1.25   \\
  1/64  & 1.9352E-02 & 1.01 & 7.5885E-02 & 1.06 & 1.8214E-01 & 1.09   \\
  1/128 & 9.6521E-02 & 1.00 & 3.7521E-02 & 1.02 & 9.0360E-02 & 1.01   \\
 \noalign{\smallskip}\hline
 \end{tabular}  \label{3dsmpnp-L2error-H1error-order-sin}
\end{table}

\subsection{Simulation on a molecular sphere model } \label{sphere-model}

\begin{figure}[H]
  \centering
  \subfigure[]{
  \label{fig:subfig:a} 
  \includegraphics[scale=0.6]{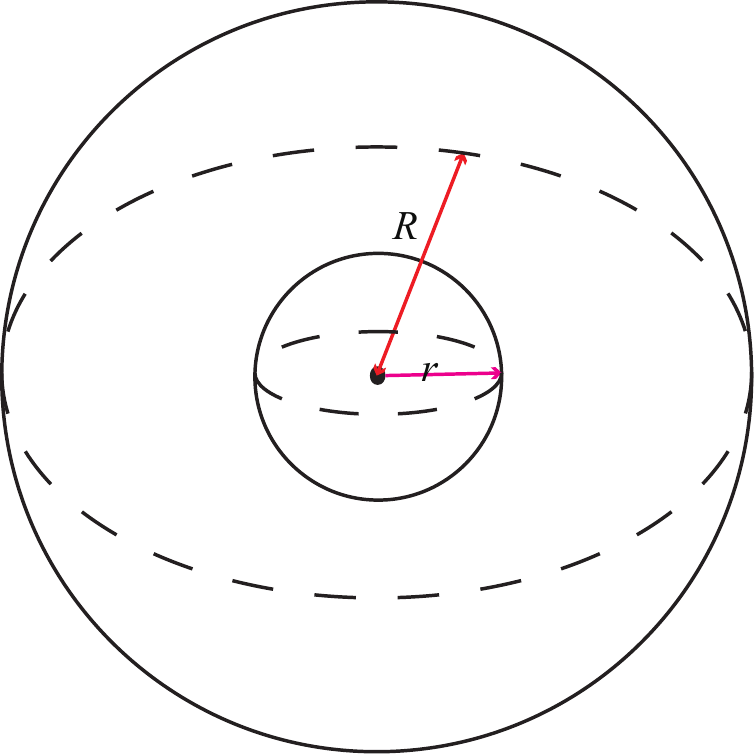} }
  \hspace{0.5in}
  \subfigure[]{
  \label{fig:subfig:b} 
  \includegraphics[scale=0.27]{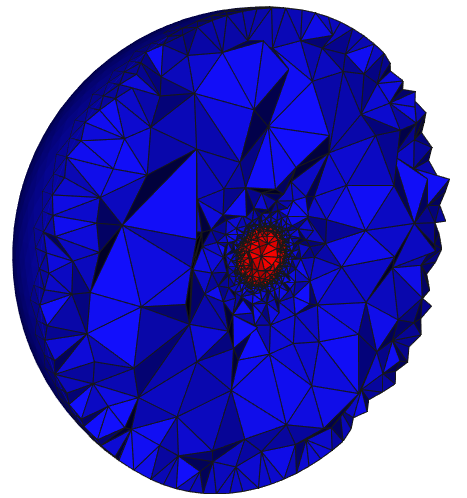} }
  \caption{ Schematic of the geometry (a) and mesh (b) of the sphere model.  }
  \label{sphere80R10r} 
\end{figure}

\noindent
In this subsection, we will use the IAFEM to solve the SMPNP/PNP equations on a sphere model. These tests capture the fundamental difference between the SMPNP model and the classical PNP model. In the sphere model, a sphere with a negative charge in the center, simulates the solute molecule. The geometry and mesh of the sphere model are shown in Fig. \ref{sphere80R10r}, where $R = 80 \AA$, $r = 10 \AA$.

$\bullet$ ~\textbf{Case 1}:

To evaluate the effectiveness and robustness of the IAFEM in bimolecule simulations, both SMPNP and PNP ($k_i = 0$, see \eqref{smpnp-nonliner}) equations are solved by setting different negative charges in the center of the molecular sphere. We first consider a $1:1$ KCl solution in our sphere model, where the bulk concentration $c_{bulk} = 0.1 M$ and the applied potential $\phi_0 = 0 V$.  The negative charges in the center of the sphere are set as: $q_1 = -10e_c$, $q_2 = -20e_c$, $q_3 = -30e_c $, $q_4 = -35e_c $, $q_5 = -40e_c $, and $q_6 = -45e_c $, where $e_c$ is the elementary charge.

\begin{figure}[H]
  \centering
  \includegraphics[scale=0.45]{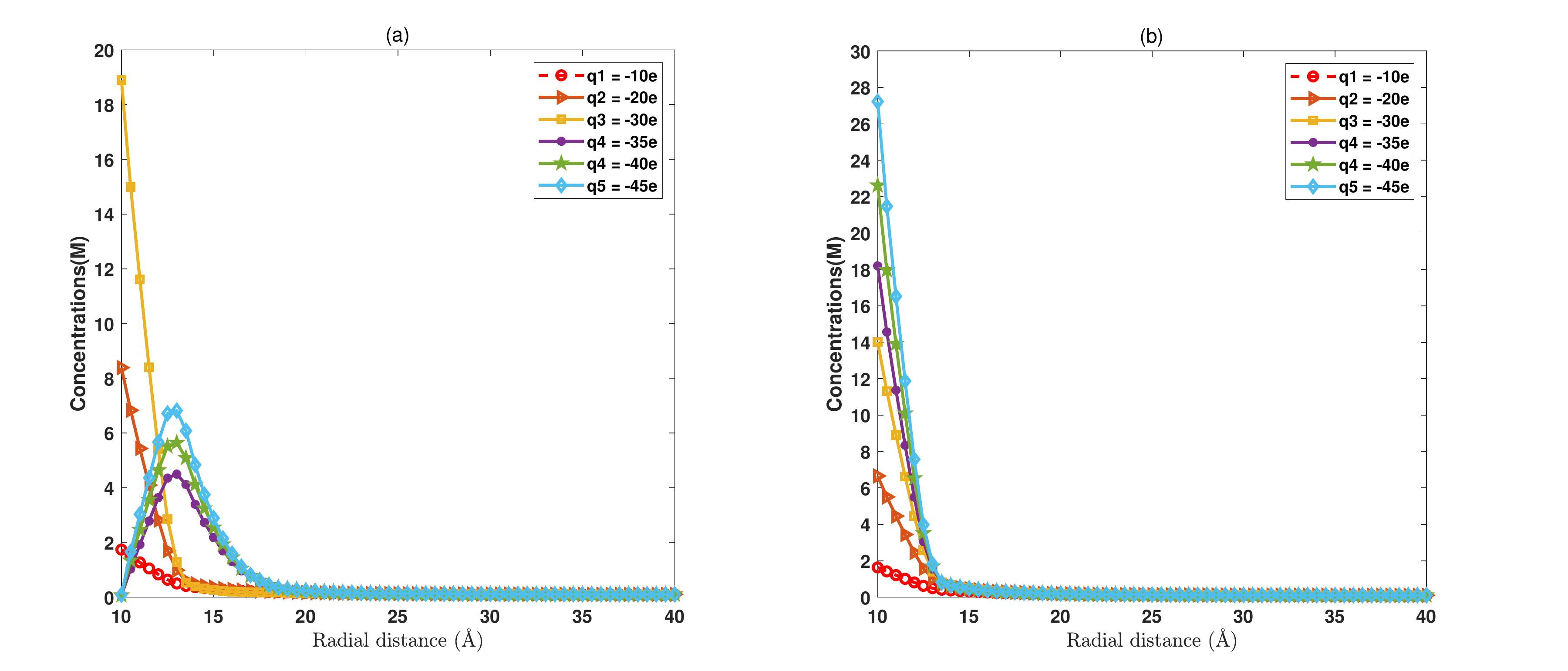}
  \caption{ Cation ($K^+$) distributions near the surface of the sphere under different charges based on PNP equations which are solved by FEM (a) and IAFEM (b). }
  \label{pnp-FE-EAFF-251K-637Cl-changeQ_31a0_01M}
\end{figure}

Our previous numerical experience shows that the traditional finite element schemes for solving PNP equations often lead to nonphysical oscillations (negative concentration values) in practical computations, e.g. see \cite{Wang2021A}. As we all know, in the simulation of the molecular spheres, the counter-ion concentration near the surface of the  molecular sphere increases as charges on the sphere accumulate. The counter-ion ($K^+$) concentrations with different charge quantity on the center sphere
solved from the classical PNP equations with the standard finite element method (FEM) and IAFEM are respectively displayed in Fig. \ref{pnp-FE-EAFF-251K-637Cl-changeQ_31a0_01M}. When using the FEM to solve the PNP equations, Fig. \ref{pnp-FE-EAFF-251K-637Cl-changeQ_31a0_01M} (a) shows that the counter-ion concentration appears layer as the amount of charges increases to a certain value, e.g., $q_4 = -35e_c $, $q_5 = -40e_c $, $q_6 = -45e_c$.
However, it's an impossible phenomenon for a single univalent ion without competition in these mean field models (the counter-ion concentration should decrease monotonically in this case).
In other words, these layers are nonphysical solutions caused by traditional FEM. Compared to Fig. \ref{pnp-FE-EAFF-251K-637Cl-changeQ_31a0_01M} (a), it is apparent from Fig. \ref{pnp-FE-EAFF-251K-637Cl-changeQ_31a0_01M} (b) that the curves of the counter-ion concentrations are always monotonous as the amount of charges increases. And the concentration value decreases rapidly to the bulk concentration value with the increase of the radial distance.

\begin{figure}[H]
  \centering
  \includegraphics[scale=0.45]{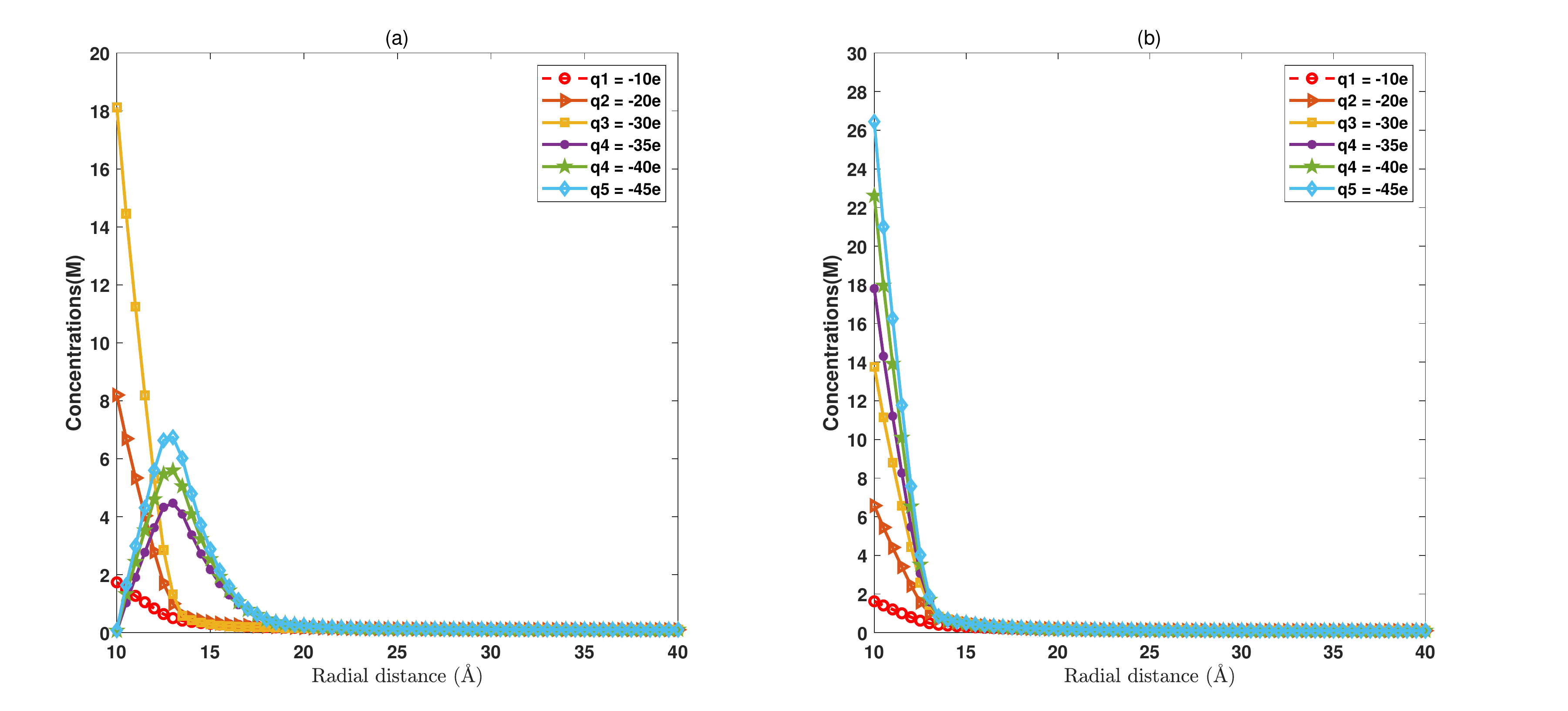}
  \caption{ Cation ($K^+$) distributions near the surface of the sphere under different charges based on SMPNP equations which are solved by FEM (a) and IAFEM (b). The solvent molecular size is $3.1 \AA$, the  counterion and coion size are $a_{K} = 2.51 \AA$ and $a_{Cl} = 6.37 \AA$, respectively. The ionic bulk densities are $0.1 M$. }
  \label{smpnp-aFE-bEAFF-251K-637Cl-changeQ_31a0_01M} 
\end{figure}

$\bullet$ ~\textbf{Case 2}:

Similarly, in order to demonstrate the effectiveness of the IAFEM for solving the SMPNP model in biomolecule simulations. In this case, we choose the ion size of cation and anion as $a_{K} = 2.51 \AA$ and $a_{Cl} = 6.37 \AA$, respectively, and the other parameters are the same as that mentioned in \textbf{Case 1}. The numerical results are shown in Fig. \ref{smpnp-aFE-bEAFF-251K-637Cl-changeQ_31a0_01M}. The curves in (a) are obtained from FEM with the standard scheme \eqref{std-smpnp-np-stdFE} for SMNP equations, and the curves in (b) are computed with IAFEM. Similarly, it is observed from Fig. \ref{smpnp-aFE-bEAFF-251K-637Cl-changeQ_31a0_01M} (a) that there also exist layers (nonphysical solutions) when solving SMPNP equations with FEM when the amount of charges increases bigger than a certain value. This illustrates that the SMPNP equations almost degenerate to PNP equations for counterion when the size of the counterion is less than the size of the solvent molecular. Fig. \ref{smpnp-aFE-bEAFF-251K-637Cl-changeQ_31a0_01M} (b) shows that if the IAFEM is used to solve the SMPNP equations, the layer will not appear even if the quantity charges is high. The effectiveness and robustness of the IAFEM are further verified for solving SMPNP/PNP equations in biomolecule simulation. In addition, comparing with Fig. \ref{pnp-FE-EAFF-251K-637Cl-changeQ_31a0_01M} (b), from Fig. \ref{smpnp-aFE-bEAFF-251K-637Cl-changeQ_31a0_01M} (b), it is seen that the counter-ion concentration solved from SMPNP equations is less than that obtained from PNP equations at the same amount of charges because of the ion size-effects.

\begin{figure}[H]
  \centering
  \includegraphics[scale=0.35]{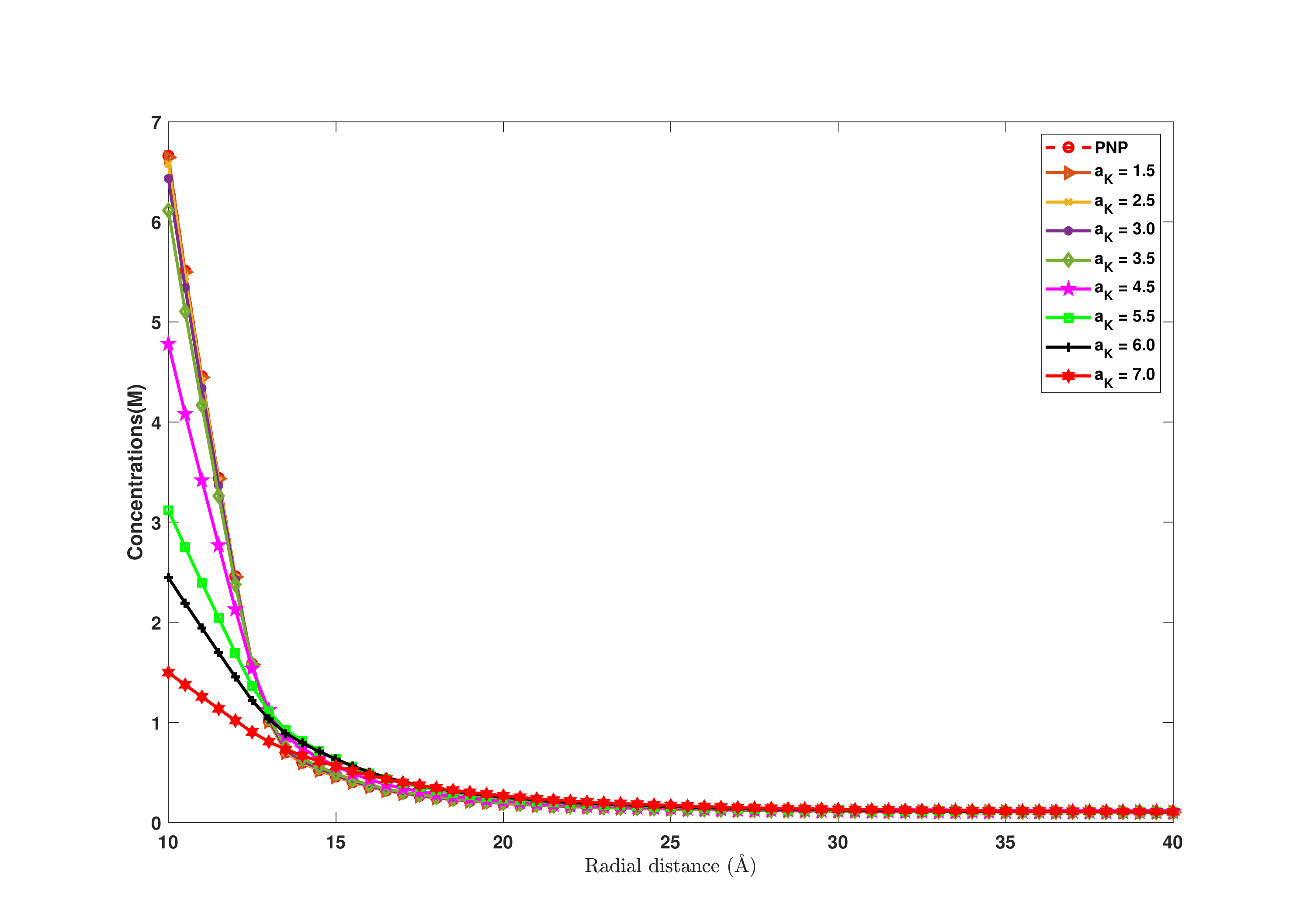}
  \caption{Cation ($K^+$) distributions near the surface of the sphere based on PNP equations (dashed line) and SMPNP equations under different counterion size (solid lines). Both the PNP and SMPNP equations are solved by IAFEM. The solvent molecular size is $3.1 \AA$, the  coion size $a_{Cl} = 6.37 \AA$. The ionic bulk densities are $0.1 M$.}
  \label{smpnp_sphere_K_637Cl_r10_R80_qf20_31a0_01M} 
\end{figure}

We know that the ion size effect has a certain inhibitory effect on the ion concentration distribution \cite{Lu2011smpnp}. In order to further reflect the influence of the ion size-effects on counterion concentration distributions, we fix the center charge of the sphere $q = - 20 e_e$ and the anion size $a_{Cl} = 6.37 \AA$. Then the concentration distributions of counterions are studied by changing the size of the counterions. The results based on the PNP equations and SMPNP equations with different counterions sizes are listed in Fig. \ref{smpnp_sphere_K_637Cl_r10_R80_qf20_31a0_01M}. These models are solved by IAFEM. When the ion size is larger, the concentration of the counterion near the surface of the sphere is smaller because of the inhibition of ion size effects.

$\bullet$ ~\textbf{Case 3}:

In this test, a $1:1:2$ mixed solution of $Na^+$, $K^+$  and $Cl^-$ is taken into account, in which the bulk concentration is ~$c_{Na^+} = c_{K^+} = 0.1 M$, $c_{Cl^-} = 0.2 M$, and the diffusion coefficients are $D_{Na^+} = 0.133 \AA/s$, $D_{K^+} = 0.196 \AA/s$, $D_{Cl^-} = 0.203 \AA/s$ for  $Na^+$, $K^+$  and $Cl^-$, respectively. For the ion size, we consider the hydration layer diameter of ions, that is $a_{Na^+} = 4.79 \AA$, $a_{K^+}=5.51\AA$, $a_{Cl^-} = 6.37 \AA$ (cf.\cite{Libo2013Ionic,Qiao2016A}). Similar to \textbf{Case 1},  we investigate the convergence of the traditional FEM and IAFEM  with various quantities of charges in the center of the sphere. In particular, in order to eliminate the influence of relevant factors of Gummel iteration on iterative convergence, we set the relaxation parameter $\alpha = 0.1$. The error tolerance $tol$ is set as $1.0 \times 10^{-6}$. If $\frac{\|u - u_{old}\|_2} {\|u\|_2} < tol$, which represents the nonlinear iterative convergence, then the computation stops.  The convergence results for FEM and IAFEM with different amounts of charges in the center of the sphere are listed in Table \ref{sphere-change-charge-3ions}, where the notations ``$\checkmark$" and ``$\times$"  represent whether the Gummel iteration has converged. The ellipsis (...) in the fifth column represents that the  Gummel iteration can converge normally when the quantity of charges is in the interval $(-25e_c, -28e_c)$ both for the traditional FEM and IAFEM. Compared with the traditional FEM, Table \ref{sphere-change-charge-3ions} indicates that the IAFEM can simulate the highly charged molecular sphere with the size effect of the hydration layer of ions for a mixed solution.

\begin{table}[H]
 \centering
 \caption{ The convergence of Gummel iteration (convergence: $\checkmark$, non-convergence: $\times$). }
 \vskip 0.1cm
 \begin{tabular}{p{1.5cm} c c c c c c c c c c c cc c c c  }
 \hline\noalign{\smallskip}
  $q$ & $-10e_c$ &$-20e_c$ & $-25e_c$&...&$-28e_c$ &$-29e_c$  & $-30e_c$ & $-35e_c$ & $-38e_c$ \\
 \noalign{\smallskip}\midrule[1pt]\noalign{\smallskip}
  FE   &  \checkmark & \checkmark & \checkmark &\checkmark & \checkmark & $\times$ & $\times$ & $\times$ & $\times$ \\
  IAFE &  \checkmark &\checkmark & \checkmark &\checkmark & \checkmark & \checkmark & \checkmark  & \checkmark & \checkmark \\
 \noalign{\smallskip}\hline
 \end{tabular}  \label{sphere-change-charge-3ions}
\end{table}

Furthermore, in order to investigate the influence of the size effect of the counterion on the robustness of our methods, we set the charge amount $q = -20e_c$ and the size of the coion $a_{Cl} = 6.37 \AA$, and the SMPNP equations are solved based on FEM and IAFEM with various sizes of the counter-ions. We choose the hydration layer diameter of $Na^+$ and $K^+$, that is $(a_{Na^+}, ~ a_{K^+}) = (4.79\AA, ~ 5.51\AA)$, as the starting point and increase the ion size by one unit ($1 \AA$) at a time. The convergence results of the traditional FEM and IAFEM are reported in Table \ref{sphere-change-aNAaK-3ions}.  It is seen from Table \ref{sphere-change-aNAaK-3ions} that the IAFEM can simulate the SMPNP equations with the strong counter-ion size effect (the counter-ion size is greater than $10\AA$) under some proper conditions. However, the traditional FEM can solve the SMPNP equations only with weak counterion size effects under the same conditions. These numerical experiments further verify the robustness and effectiveness of the IAFEM for solving the SMPNP equations with large-size effects in biomolecule simulations.

\begin{table}[H]
 \centering
 \caption{ The convergence of Gummel iteration (convergence: $\checkmark$, non-convergence: $\times$) }
 \vskip 0.1cm
 \begin{tabular}{p{2.5cm} c c c c c   }
 \hline\noalign{\smallskip}
  $(a_{Na^+}, a_{K^+})(\AA)$ & $(4.79, 5.51)$ & $(5.79, 6.51)$ & $(6.79, 7.51)$ &$(7.79, 8.51)$ & $(8.79, 9.51)$\\
 \noalign{\smallskip}\midrule[1pt]\noalign{\smallskip}
  FE   & \checkmark & \checkmark  & \checkmark & \checkmark & $\times$   \\
  IAFE & \checkmark & \checkmark  & \checkmark & \checkmark & \checkmark  \\
  \hline
 \hline\noalign{\smallskip}
  $(a_{Na^+}, a_{K^+})(\AA)$ &$(9.79, 10.51)$ & $(10.79, 11.51)$ & $(11.79, 12.51)$ & $(12.79, 13.51)$ \\
 \noalign{\smallskip}\midrule[1pt]\noalign{\smallskip}
  FE   & $\times$   & $\times$   & $\times$   & $\times$ \\
  IAFE & \checkmark & \checkmark & \checkmark & \checkmark \\
 \noalign{\smallskip}\hline
 \end{tabular}  \label{sphere-change-aNAaK-3ions}
\end{table}

\subsection{Size-effects in ion transports: a numerical simulation of a gA channel  }

\begin{figure}[H]
 \centerline{
 \includegraphics[scale=0.3]{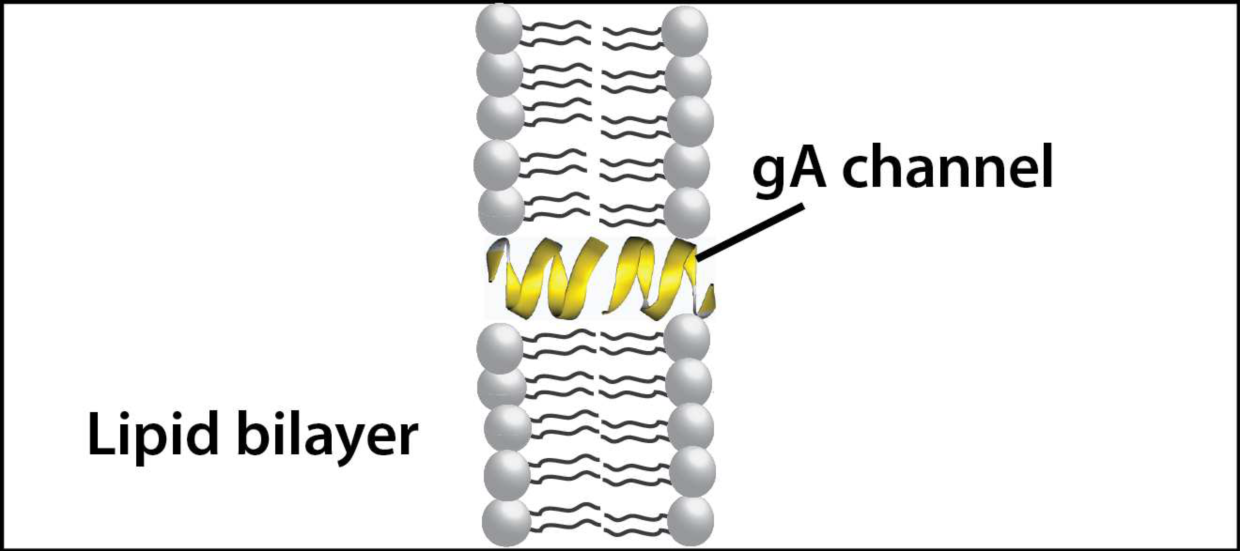} }
 \caption{A schematic picture of a single gramicidin A pore embedded in a lipid bilayer (see \cite{Xu2014Modeling}). }
 \label{gA-lipid-bilayer}
\end{figure}

\noindent
 In this subsection, we use the IAFEM to solve the SMPNP equations and PNP equations in the simulation of an ion channel. Gramicidin A (gA) is a well-characterized short polypeptide including hundreds of atoms with a helix structure.  Fig. \ref{gA-lipid-bilayer} shows a schematic picture of a single gA channel embedded in the lipid bilayer. This peptide is relatively easy to be synthesized and manipulated, compared with a typical sodium channel which has thousands of atoms. The gA channel is also relatively stable, therefore, it has been widely applied in biochemical and biophysical studies. Upon head to head dimerization, gA forms an elongated channel in the lipid bilayer that is permeable to small monovalent cations \cite{Wallace1986Structure}.

\begin{figure}[H]
 \centerline{
 \includegraphics[scale=0.81]{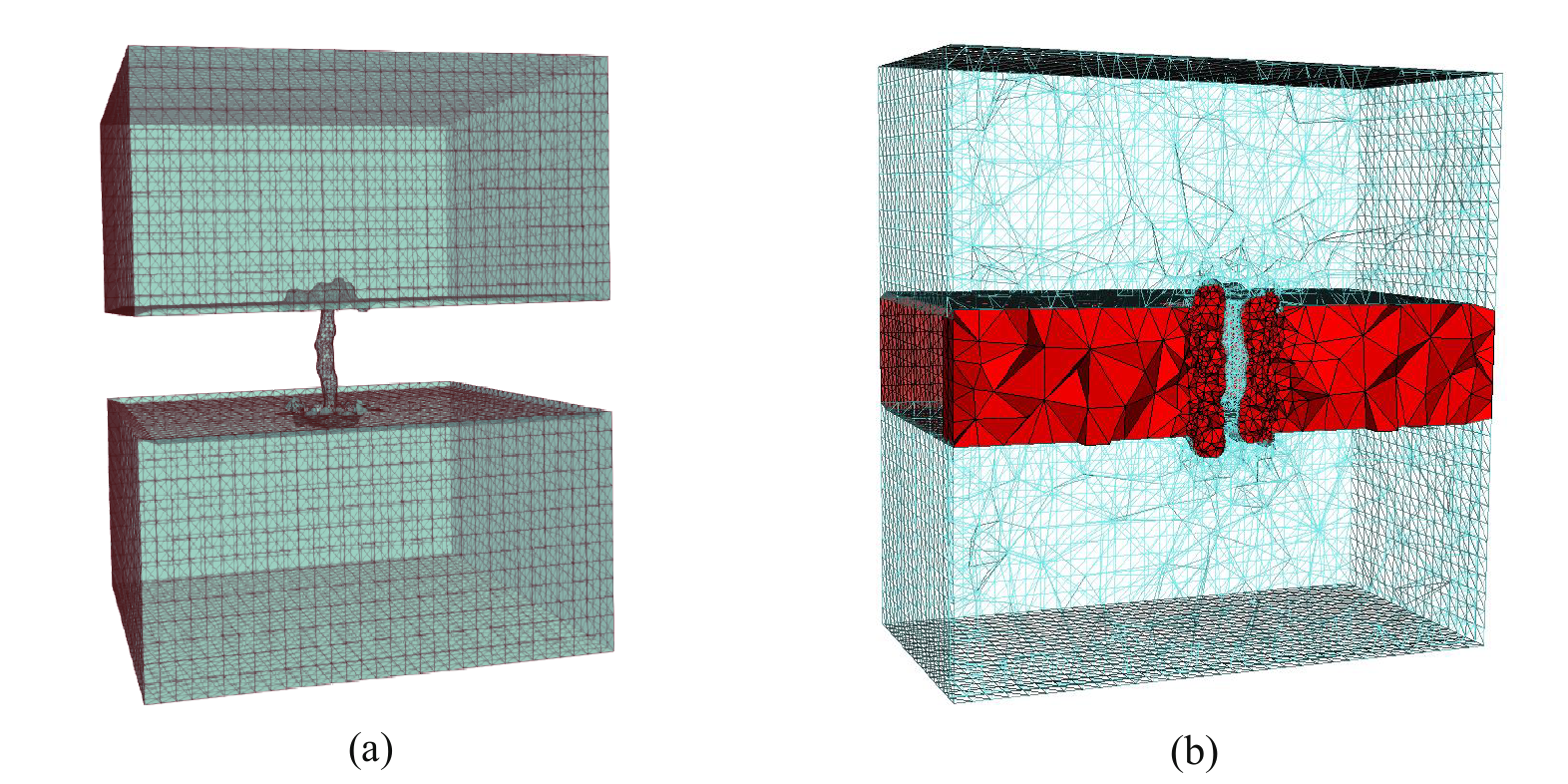} }
 \caption{ a) The gA channel embedded in the simulation box. b) A cut plane through the center of the simulation box along the z-axis.}
 \label{gA-box-mesh}
\end{figure}

\begin{figure}[H]
 \centerline{
 \includegraphics[scale=0.5]{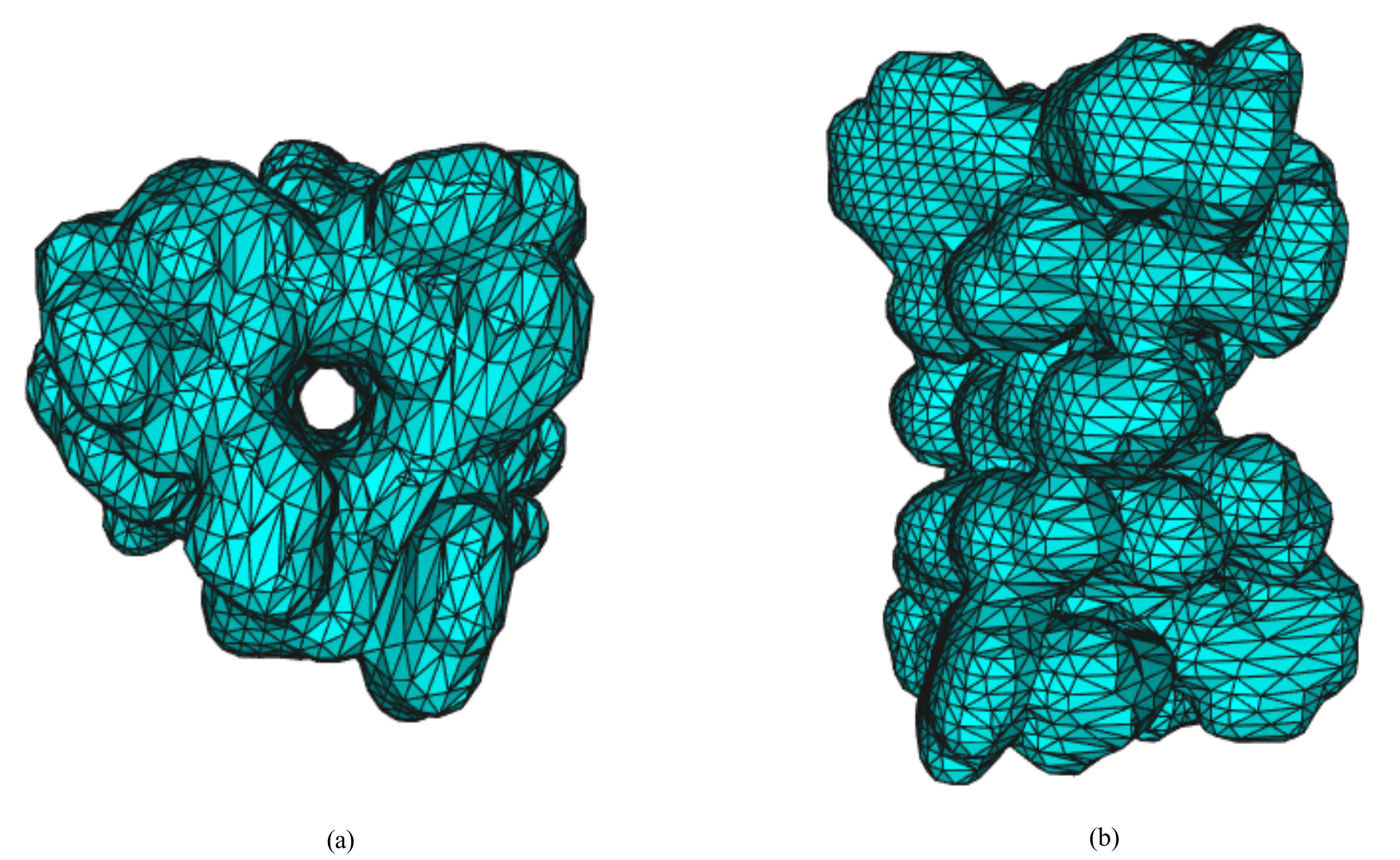} }
 \caption{ Schematic of the gA channel meshes:  a) top view. b) lateral view.}
 \label{gA-meshes}
\end{figure}

 In our work, we utilize the SMPNP and PNP equations to calculate the concentration distributions of the cations in the channel. The size effects on the ion concentration distributions in the channel are further studied for different ions. In our computation, the gA channel system setup is similar to the model presented in \cite{Tu2013A}. The whole domain of the gA channel consists of the membrane protein region, bulk region, and channel region. The gA channel region is along the $z$-direction. The schematic of the simulation box and the gA channel meshes are shown in Figs. \ref{gA-box-mesh} and \ref{gA-meshes}. In our simulations, the box size is $[-50\AA, 50\AA]^3$, the channel region is $[-14\AA,7\AA]$, and the total number of tetrahedral elements is 92480.

\begin{figure}[H]
 \centerline{
 \includegraphics[scale=0.45]{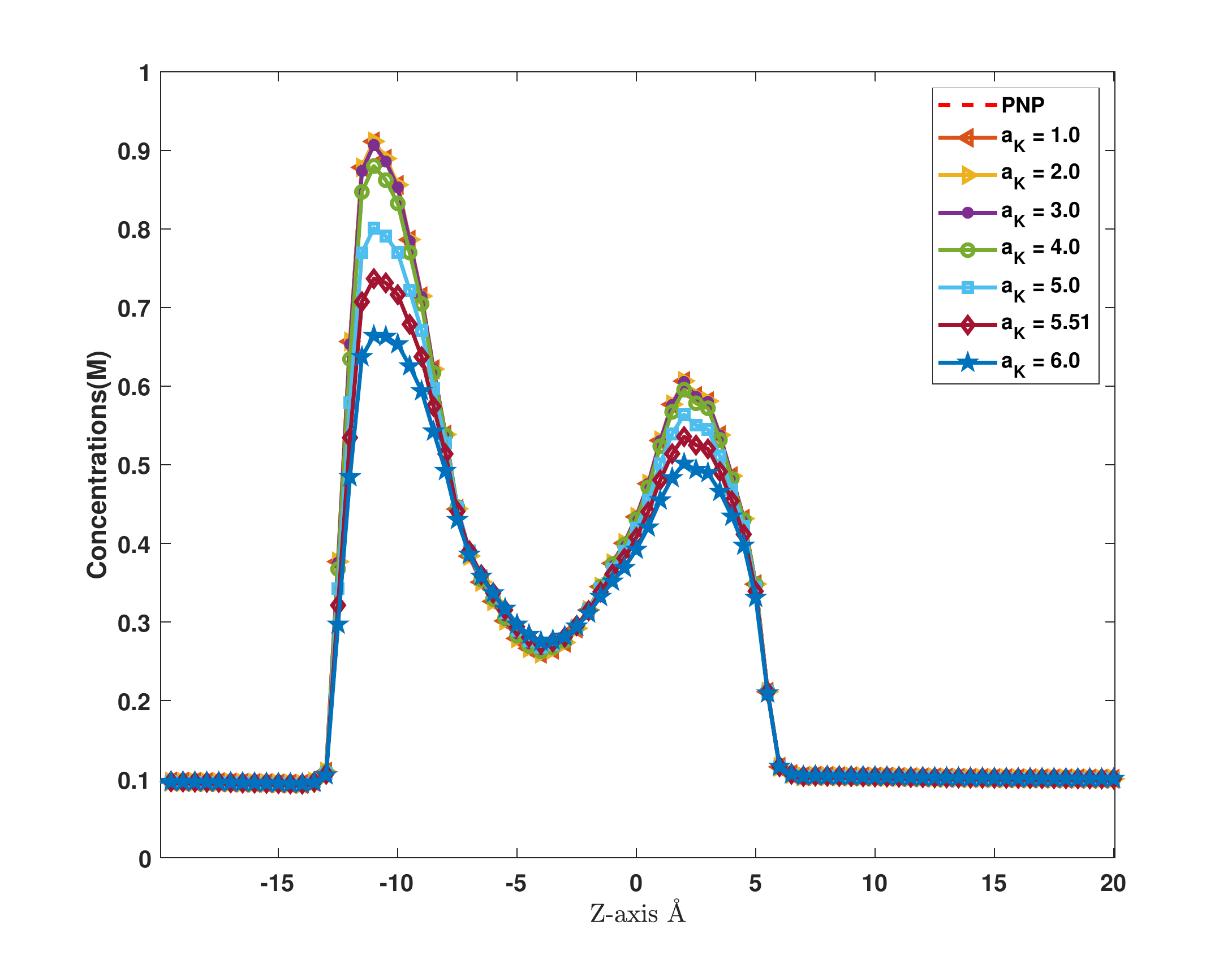} }
 \caption{ The cation density profile in the channel computed by PNP equations (dashed line) and SMPNP equations with different cation size (solid lines), where the anion size $a_{Cl}= 4.37 \AA$, and the solvent molecular size $a_0 = 3.1 \AA$.  }
 \label{smpnp-K-gA-con01-minus-f015}
\end{figure}

$\bullet$ ~\textbf{Case 1}:

At first, similar to the sphere model, we consider a $1:1$ KCl solution in our gA channel system, where the bulk concentration $c_{bulk} = 0.1 M$ is on the top and bottom of the box, and the potential $\phi_0 = -0.15V$ is applied with the potential difference along the z-direction. The diffusion coefficients for cation and anion, for example, $K^{+}$ and $Cl^{-}$, in the bulk region are set to their experimental values: $D_{K}=0.196 \AA^{2} / ps$, $D_{\mathrm{Cl}}=0.203 \AA^{2} / \mathrm{ps}$. While there is no experimental measurement of exact values for the diffusion coefficients inside the channel, it is known that the diffusion coefficients in the bulk region and the channel region should be different. In this work, the diffusion coefficients inside the channel are set by the same way as shown in \cite{Tu2013A}. In order to study the influence of size effects on the cation concentration distribution in the channel, both the PNP equations and SMPNP equations with different ionic sizes are solved by the IAFEM for the gA system. The size of the anion ($Cl^-$) is fixed at  $4.37 \AA$, and the solvent molecular size $a_0 = 3.1 \AA$. The cation size is arbitrarily given in this test. Fig. \ref{smpnp-K-gA-con01-minus-f015} shows the cation concentration profile inside the channel  solved from PNP equations and SMPNP equations with different cation sizes. If the size of the cation is larger, the cation concentration in the channel obtained from SMPNP equations is lower compared with that solved from PNP equations. This is because the large size of the cation can strengthen the spatial repulsion of the model. Therefore, the SMPNP model, especially with large size effects, can be used to control the infinite increase of ion concentrations in the channel. For example, the ion currents matched with the experimental data (cf. \cite{Andersen1983Ion-1,Andersen1983Ion-2,Andersen1983Ion-3}) can be obtained through numerical simulations with high bulk concentrations and high applied voltage difference, which will be studied in our next work.

$\bullet$ ~\textbf{Case 2}:

 In this case, we consider the influence of the size effects on ion competitions in the channel. Similar to the test in \textbf{Case 3} of the biomolecular sphere model mentioned above, we still consider a $1:1:2$ mixed salt solution for $Na^+$, $K^+$  and $Cl^-$, in which the bulk concentration is ~$c_{Na^+} = c_{K^+} = 0.1 M$, $c_{Cl^-} = 0.2 M$, and the diffusion coefficients and the applied potential are the same as that in \textbf{Case 1}. In addition, the solvent molecular size and the anion size are respectively fixed as $a_0 = 3.1 \AA$ and $a_{Cl}= 4.37 \AA$. Both the PNP and SMPNP equations are solved via the IAFEM.

\begin{figure}[H]
\centerline{
 \includegraphics[scale=0.45]{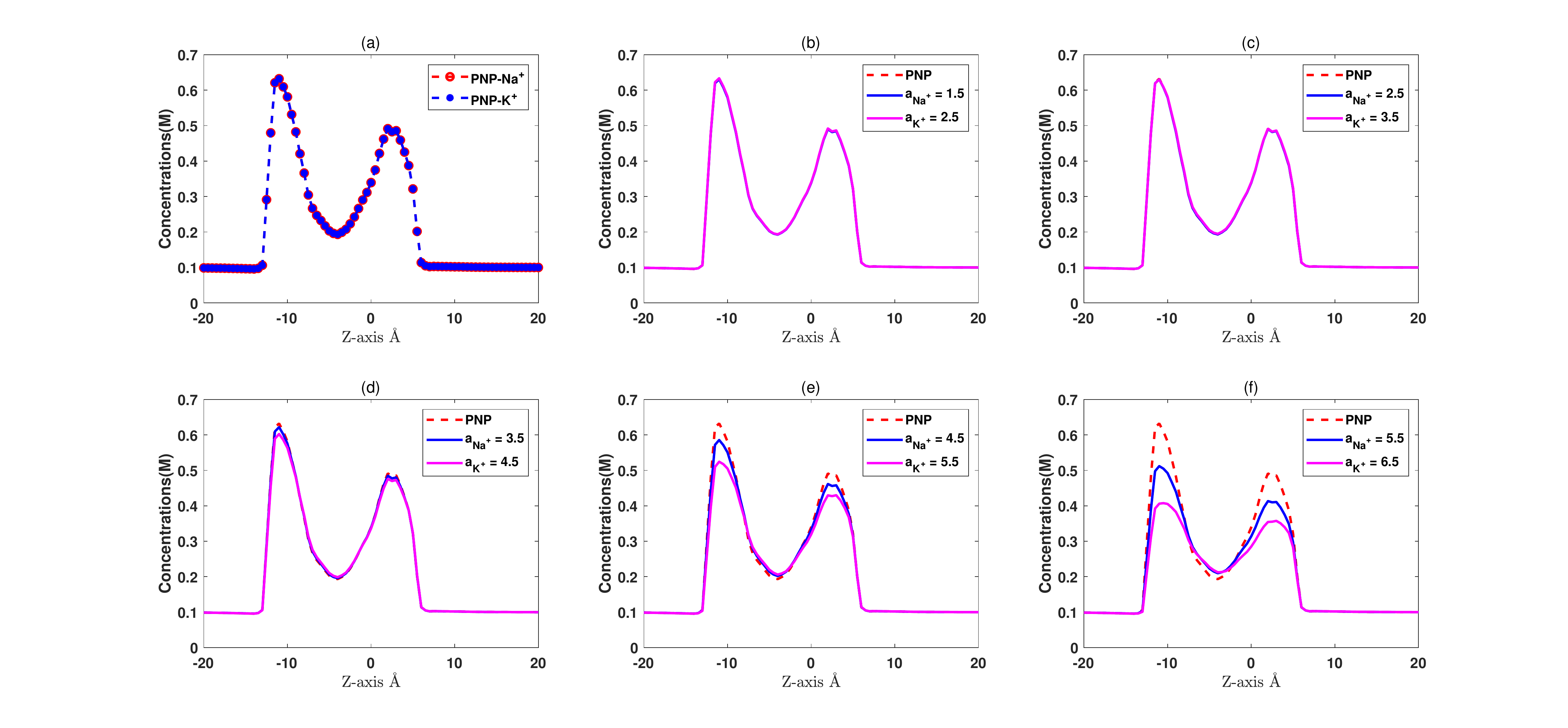} }
 \caption{Cation distributions in the gA channel under a fixed membrane voltage ($\phi_0 = -0.15V $) and bulk concentration ($c_{Na^+} = c_{K^+} = 0.1 M$, $c_{Cl^-} = 0.2 M$) computed by PNP equations (dashed line) and SMPNP equations with different ion size (solid lines). (a): The ion distributions both for $Na^+$ and $K^+$ by PNP equations; (b) - (f): The ion distributions based on SMPNP equations with different cation size (solid lines). The anion size $a_{Cl}$ and the solvent molecular size $a_0$ is fixed as $4.37 \AA$ and $3.1\AA$, respectively. }
 \label{smpnp-K-Na-gA-size-con01-minus-f015}
\end{figure}

 For SMPNP equations, to evaluate the influence of ion size effects on cation distributions in the channel, we arbitrarily change and increase the sizes of $Na^+$ and $K^+$ in our computation. It is worth noting that the actual size of $K^+$ is larger than the size of $Na^+$, so we always make the ion size for $K^+$ one unit ($1 \AA$) larger than the size for $Na^+$ at a time in the test. We considered the concentration distributions in the gA channel along the $z$-axis. The numerical results are displayed in Fig. \ref{smpnp-K-Na-gA-size-con01-minus-f015}. The subfigure (a) is the concentration profile for $Na^+$ and $K^+$ obtained by PNP equations without the size effects, and the subfigures (b) - (f) show the concentration profile obtained from SMPNPEs with different ion sizes for $Na^+$ and $K^+$. Fig. \ref{smpnp-K-Na-gA-size-con01-minus-f015} (a) demonstrates that the traditional PNP equations cannot distinguish the two cations with the same valence. In particular, for comparison, the concentration distributions computed through PNP equations are also shown in figures (b) - (f). It is observed from Fig. \ref{smpnp-K-Na-gA-size-con01-minus-f015} that, similar to PNP equations, the ion distributions for $Na^+$ and $K^+$ calculated with SMPNP equations still cannot be distinguished in the channel when the ion size is small (see subfigures (b) - (d)). However, the subfigures (e) - (f) show that the concentration distributions of $Na^+$ and $K^+$ in the channel can be clearly distinguished by SMPNP equations when the cation size becomes larger. In addition, the concentration profiles from subfigures (e) - (f) also show that the concentration of $K^+$ is less than that of $Na^+$, which indicates that the size of the ion can strengthen the size inhibition effect.

%
%
%
%

\section{Conclusion} \label{sec-conclusion}

\noindent
 In this work, we introduce a generalized Slotboom transformation and an IAFEM to solve the SMPNP equations. With the generalized Slotboom transform, the original SMNP equations are transformed into new reformulations which are self-adjoint equations with exponentially behaved coefficients. Then the inverse averaging technique over the edges of the element can be used to deal with  the exponential coefficients.  Numerical experiments are reported to confirm the efficiency and robustness of the new schemes for SMPNP equations. Firstly, a model problem with analytic solutions on a cube box is tested to numerically verify the accuracy and order of the new schemes for SMPNP equations. Then, based on the averaging technique, simulations of both biomolecular sphere systems and ion channel systems are studied to demonstrate the effectiveness and robustness of the IAFEM for SMPNP equations. In particular, the averaging technique introduced in this paper can be easily extended to more complex PNP-like models for simulating biomolecular systems, such as the variable dielectric Poisson-Nernst-Planck (VDPNP) equations and Born-energy-modified PNP (BPNP) equations \cite{Liu2017Incorporating,Liu2018Analysis}, in which the dielectric coefficients depend on the ion concentrations and spatial positions, respectively. This will be studied in our future work. In addition, because of the strong coupling non-linearity of the system and the proposed scheme is based on the reformulation of the SMNP equations with exponential terms, the convergence analysis is not a matter of standard analysis. We leave it as the future work.

\section*{Acknowledgments}
Many thanks for the help and discussions with Dr. Y. Qiao. Thanks also go to Dr. S. Gui for his help in mesh generation and visualization.  B. Z. Lu was supported by the China NSF (NSFC 22073110, NSFC 11771435). R. G. Shen was supported by the China NSF (NSFC 12101595).




\begin{thebibliography}{99} \setlength{\itemsep}{- 0.5mm}
\bibitem{Andersen1983Ion-1}
O.~S. Andersen.
 Ion movement through gramicidin A channels. Single-channel
  measurements at very high potentials.
 {\em Biophys. J.}, 41(2):119--133, 1983.

\bibitem{Andersen1983Ion-2}
O.~S. Andersen.
 Ion movement through gramicidin A channels. Interfacial polarization
  effects on single-channel current measurements.
 {\em Biophys. J.}, 41(2):135--146, 1983.

\bibitem{Andersen1983Ion-3}
O.~S. Andersen.
 Ion movement through gramicidin A channels. Studies on the
  diffusion-controlled association step.
 {\em Biophys. J.}, 41(2):147--165, 1983.

\bibitem{Bazant2009Towards}
M.~Z. Bazant, M.~S. Kilic, B.~D. Storey, and A.~Ajdari.
 Towards an understanding of induced-charge electrokinetics at large
  applied voltages in concentrated solutions.
 {\em Adv. Colloid Interfac.}, 152(1-2):48--88, 2009.

\bibitem{Borukhov1998Steric}
I.~Borukhov, D.~Andelman, and H.~Orland.
 Steric effects in electrolytes: A modified Poisson-Boltzmann
  equation.
 {\em Phys. Rev. Lett.}, 79(3):435--438, 1998.

\bibitem{Burman2004Edge}
E.~Burman and P.~Hansbo.
 Edge stabilization for galerkin approximations of
  convection-diffusion-reaction problems.
 {\em Comput. Methods Appl. Mech. Engrg.}, 193(15-16):1437--1453,
  2004.

\bibitem{AE2000Three}
A.~E. C{\'a}rdenas, R.~D. Coalson, and M.~G. Kurnikova.
 Three-dimensional Poisson-Nernst-Planck theory studies: influence of
  membrane electrostatics on gramicidin A channel conductance.
 {\em Biophys. J.}, 79(1):80--93, 2000.

\bibitem{Chaudhry2014A}
J.~H. Chaudhry, J.~Comer, A.~Aksimentiev, and L.~N. Olson.
 A stabilized finite element method for modified Poisson-Nernst-Planck
  equations to determine ion flow through a nanopore.
 {\em Commun. Comput. Phys.}, 15(1):93--125, 2014.

\bibitem{Chern2003Accurate}
I.~Liang Chern, J.~G. Liu, and W.~C. Wang.
 Accurate evaluation of electrostatics for macromolecules in solution.
 {\em Methods Appl. Anal.}, 10(2):309--328, 2003.

\bibitem{Ciucci2011Derivation}
F.~Ciucci and W.~Lai.
 Derivation of micro/macro lithium battery models from homogenization.
 {\em Transp. Porous Med.}, 88(2):249--270, 2011.

\bibitem{Coalson2005Poisson}
D.~R. Coalson and G.~M. Kurnikova.
 Poisson-Nernst-Planck theory approach to the calculation of current
  through biological ion channels.
 {\em IEEE Trans. NanoBiosci.}, 4(1):81--93, 2005.

\bibitem{Ding2019Positivity}
J.~Ding, Z. Wang, and S. Zhou.
 Positivity preserving finite difference methods for
  Poisson-Nernst-Planck equations with steric interactions: Application to
  slit-shaped nanopore conductance.
 {\em J. Comput. Phys.}, 397:108864, 2019.

\bibitem{Douglas1976Interior}
J.~Douglas and T.~Dupont.
 {\em Interior Penalty Procedures for Elliptic and Parabolic Galerkin Methods}.
 Springer Berlin Heidelberg, 1976.

\bibitem{Flavell2014A}
A. Flavell, M. Machen, B. Eisenberg, J. Kabre, C. Liu, and X. Li.
 A conservative finite difference scheme for Poisson-Nernst-Planck equations.
 {\em J. Comput. Electron.}, 13:235--249, 2014.

\bibitem{Gillespie2002Coupling}
D.~Gillespie, W.~Nonner, and R.~S. Eisenberg.
 Coupling Poisson-Nernst-Nlanck and density functional theory to calculate ion flux.
 {\em J. Phys. Condens. Mat.}, 14(46):12129--12145, 2002.

\bibitem{Gillespie2003Density}
D.~Gillespie, W.~Nonner, and R.~S. Eisenberg.
 Density functional theory of charged, hard-sphere fluids.
 {\em Phys. Rev. E Stat. Nonlin. Soft Mat. Phys.},
  68(3):031503(1--10), 2003.

\bibitem{Gummel1964A}
H.~K. Gummel.
 A self-consistent iterative scheme for one-dimensional steady state
  transistor calculations.
 {\em IEEE Trans. Electron Dev.}, 11(10):455--465, 1964.

\bibitem{He2016An}
D. He and K. Pan.
 An energy preserving finite difference scheme for the
  Poisson-Nernst-Planck system.
 {\em Appl. Math. Comput.}, 287-288:214--223, 2016.

\bibitem{Holst2012Adaptive}
M.~Holst, J.~A. Mccammon, Z.~Yu, Y.~C. Zhou, and Y.~Zhu.
 Adaptive finite element modeling techniques for the Poisson-Boltzmann equation.
 {\em Commun. Comput. Phys.}, 11(01):179--214, 2012.

\bibitem{Horng2012PNP}
T.~L. Horng, T.~C. Lin, C.~Liu, and B.~Eisenberg.
 PNP equations with steric effects: a model of ion flow through channels.
 {\em J. Phys. Chem. B}, 116(37):11422--11441, 2012.

\bibitem{Jerome1996}
J.~Jerome.
 {\em Analysis of Charge Transport: A Mathematical Theory and
  Approximation of Semiconductor Models}.
 Springer-Verlag, New York, 1996.

\bibitem{Kilic2007Steric}
M.~S. Kilic, M.~Z. Bazant, and A.~Ajdari.
 Steric effects in the dynamics of electrolytes at large applied
 voltages: II. Modified Poisson-Nernst-Planck equations.
 {\em Phys. Rev. E}, 75(2):021503, 2007.

\bibitem{Li2009Continuum}
B.~Li.
 Continuum electrostatics for ionic solutions with non-uniform ionic sizes.
 {\em Nonlinearity}, 22(4):811--833, 2009.

\bibitem{Libo2013Ionic}
B.~Li, P. Liu, Z. Xu, and S. Zhou.
 Ionic size effects: generalized Boltzmann distributions, counterion
  stratification and modified Debye length.
 {\em Nonlinearity}, 26(10):2899--2922, 2013.

\bibitem{Lin2014A}
T.~C. Lin and B. Eisenberg.
 A new approach to the lennard-jones potential and a new model:
  PNP-steric equations.
 {\em Commun. Math. Sci.}, 12(1):149--173, 2014.

\bibitem{Liu2014A}
H. Liu and Z. Wang.
 A free energy satisfying finite difference method for
  Poisson-Nernst-Planck equations.
 {\em J. Comput. Phys.}, 268(2):363--376, 2014.

\bibitem{Liu2017Incorporating}
X. Liu and B. Lu.
 Incorporating born solvation energy into the three-dimensional
  Poisson-Nernst-Planck model to study ion selectivity in KcsA ${\rm K^{+}}$ channels.
 {\em Phys. Rev. E}, 96(6):062416, 2017.

\bibitem{Liu2018Analysis}
X. Liu, Y.~Qiao, and B. Lu.
 Analysis of the mean field free energy functional of electrolyte
  solution with non-zero boundary conditions and the generalized PB/PNP
  equations with inhomogeneous dielectric permittivity.
 {\em SIAM J. Appl. Math.}, 78(2):1131--1154, 2018.

\bibitem{Lu2010pnp}
B. Lu, M.~J. Holst, J.~A. Mccammon, and Y.~C. Zhou.
 Poisson-Nernst-Planck equations for simulating biomolecular
  diffusion-reaction processes I: Finite element solutions.
 {\em J. Comput. Phys.}, 229(19):6979--6994, 2010.

\bibitem{Lu2011smpnp}
B. Lu and Y.~C. Zhou.
 Poisson-Nernst-Planck equations for simulating biomolecular
  diffusion-reaction processes II: Size effects on ionic distributions and
  diffusion-reaction rates.
 {\em Biophys. J.}, 100(10):2475--2485, 2011.

\bibitem{Lu2007Electrodiffusion}
B. Lu, Y.~C. Zhou, G.~A. Huber, S.~D. Bond, and J.~A. McCammon.
 Electrodiffusion: A continuum modeling framework for biomolecular
  systems with realistic spatiotemporal resolution.
 {\em J. Chem. Phys.}, 127(13):10B604--78, 2007.

\bibitem{Marcicki2012Comparison}
J. Marcicki, A.~T. Conlisk, and G. Rizzoni.
 Comparison of limiting descriptions of the electrical double layer
  using a simplified lithium-ion battery model.
 {\em ECS Transactions}, 41(14):9--21, 2012.

\bibitem{Markowich1986}
P.~Markowich.
 {\em The Stationary Semiconductor Device Equation}.
 Springer-Verlag, New York, 1986.

\bibitem{Outhwaite1980Theory}
C.~W. Outhwaite, L.~B. Bhuiyan, and S.~Levine.
 Theory of the electric double layer using a modified Poisson-Boltzman equation.
 {\em J. Chem. Soc., Faraday Trans. II}, 76:1388--1408, 1980.

\bibitem{Qiao2016A}
Y.~Qiao, X. Liu, M. Chen, and B. Lu.
 A local approximation of fundamental measure theory incorporated into
  three dimensional Poisson-Nernst-Planck equations to account for hard sphere
  repulsion among ions.
 {\em J. Stat. Phys.}, 163(1):156--174, 2016.

\bibitem{Rosenfeld1997Fundamental}
Y.~Rosenfeld, M.~Schmidt, H.~L{\"o}wen, and P.~Tarazona.
 Fundamental-measure free-energy density functional for hard spheres:
  Dimensional crossover and freezing.
 {\em Phys. Rev. E Stat. Phys. Plasmas Fluids Relat. Interdiscip.
  Topics}, 55(4):4245--4263, 1997.

\bibitem{Selberherr1984}
S.~Selberherr.
 {\em Analysis and Simulation of Semiconductor Devices}.
 Springer-Verlag, Wien, New York, 1984.

\bibitem{Siddiqua2017A}
F. Siddiqua, Z. Wang, and S. Zhou.
 A modified Poisson-Nernst-Planck model with excluded volume effect:
  Theory and numerical implementation.
 {\em Commun. Math. Sci.}, 16(1):251--271, 2018.

\bibitem{Singer2008A}
A.~Singer and J.~Norbury.
 A Poisson-Nernst-Planck model for biological ion channels--an
  asymptotic analysis in a three-dimensional narrow funnel.
 {\em SIAM J. Appl. Math.}, 70(3):949--968, 2009.

\bibitem{Slotboom1973Computer}
J.~W. Slotboom.
 Computer-aided two-dimensional analysis of bipolar transistors.
 {\em IEEE Trans. Electron. Devices}, 20(8):669--679, 1973.

\bibitem{Song2004Continuum}
Y. Song, Y. Zhang, C.~L. Bajaj, and N.~A. Baker.
 Continuum diffusion reaction rate calculations of wild-type and
  mutant mouse acetylcholinesterase: Adaptive finite element analysis.
 {\em Biophys. J.}, 87(3):1558--1566, 2004.

\bibitem{Tu2013A}
B. Tu, M. Chen, Y. Xie, L. Zhang, B. Eisenberg, and B. Lu.
 A parallel finite element simulator for ion transport through
  three-dimensional ion channel systems.
 {\em J. Comput. Chem.}, 34(24):2065--2078, 2013.

\bibitem{Tu2015Stabilized}
B. Tu, Y. Xie, L. Zhang, and B. Lu.
 Stabilized finite element methods to simulate the conductances of ion channels.
 {\em Comput. Phys. Commun.}, 188:131--139, 2015.

\bibitem{Wallace1986Structure}
B.~A. Wallace.
 Structure of gramicidin A.
 {\em Biophys. J.}, 49(1):295--306, 1986.

\bibitem{Wang2021A}
Q. Wang, H. Li, L. Zhang, and B. Lu.
 A stabilized finite element method for the Poisson-Nernst-Planck
  equations in three-dimensional ion channel simulations.
 {\em Appl. Math. Lett.}, 111:106652(1--9), 2021.

\bibitem{Xu2014Modeling}
S. Xu, M. Chen, S. Majd, X. Yue, and C. Liu.
 Modeling and simulating asymmetrical conductance changes in
  gramicidin pores.
 {\em Molecular Based Mathematical Biology}, 2(1):34--55, 2014.

\bibitem{Zhang2009A}
L. Zhang.
 A parallel algorithm for adaptive local refinement of tetrahedral
  meshes using bisection.
 {\em Numer. Math. Theor. Meth. Appl.}, 02(1):65--89, 2009.

\bibitem{ZhangQianru2021A}
Q. Zhanng, Q. Wang, B. Lu, and L. Zhang.
 A class of finite element methods with averaging techniques for
  solving the three-dimensional drift-diffusion model in semiconductor device
  simulations.
 {\em J. Comput. Phys.}, submitted:to appear, 2021.

\bibitem{Zhanng2021An}
Q. Zhanng, Q. Wang, L. Zhang, and B. Lu.
 An inverse averaging finite element method for solving
  three-dimensional Poisson-Nernst-Planck equations in nanopore system
  simulations.
 {\em J. Chem. Phys.}, 155(19):194106, 2021.



\end{thebibliography}



\end{document}